\documentclass{amsart}
\usepackage{graphicx} 
\usepackage{aaPreamble}
\usepackage{xcolor}
\usepackage{mathtools,amsmath,amsthm,amsfonts,amssymb,fancyhdr,graphics}
\usepackage{comment}
\usepackage{hyperref}
\usepackage{tikz-cd}
\usepackage{array}
\usepackage[all]{xy}
\usepackage{longtable}

\def\tiago#1{{\color{blue}\textbf{Tiago: }#1}}

\numberwithin{equation}{section}

\theoremstyle{plain}
\newtheorem{lemma}{Lemma}[section]

\newtheorem{theorem}[lemma]{Theorem}

\theoremstyle{definition}

\newtheorem*{definition*}{Definition}
\newtheorem*{remark*}{Remark}
\newtheorem*{remarks*}{Remarks}
\newtheorem*{example*}{Example}
\newtheorem*{examples*}{Examples}
\newtheorem*{conjecture*}{Conjecture}
\newtheorem*{conjectures*}{Conjectures}
\newtheorem*{exercise*}{Exercise}
\newtheorem*{exercises*}{Exercises}
\newtheorem*{problem*}{Problem}
\newtheorem*{problems*}{Problems}

\title{A survey on birational Rigidity of threefold Weighted Complete Intersections}

\author{Tiago Duarte Guerreiro}
\address[Tiago Duarte Guerreiro]{Departement Mathematik und Informatik, Universität Basel, Spiegelgasse 1, 4051 Basel, Switzerland}
\email{tiago.duarteguerreiro@unibas.ch}

\author{Takuzo Okada}
\address[Takuzo Okada]{Faculty of Mathematics, Kyushu University, Fukuoka 819-0385, Japan}
\email{tokada@math.kyushu-u.ac.jp}
\date{\today}
\begin{document}
\maketitle

\begin{center}
\emph{To Miles Reid, \\ for showing us the way.}
\end{center}

\begin{abstract}
We survey what is known about Fano threefold weighted complete intersections from the point of view of birational rigidity.

\end{abstract}

\tableofcontents

\section{Introduction}

One of the biggest achievements in birational geometry has been the possibility of extension of the birational classification of algebraic surfaces to threefolds and beyond, as envisioned by Shigefumi Mori in the seventies, via the Minimal Model Program (MMP). Very roughly, the MMP aims at constructing simpler representatives in the birational class of a given variety. 

A major step forward was achieved in \cite{BCHM}, where, in particular, it has been proven that MMP holds for uniruled varieties, that is, varieties covered by rational curves. In this case, a representative of the class is a \emph{Mori fibre space} which can be thought of as a relative generalisation of a Fano variety. More concretely, a Mori fibre space is an extremal contraction $\varphi \colon Y \rightarrow B$ with connected fibres from a $\mathbb{Q}$-factorial terminal projective variety $Y$ to a normal projective variety $B$ whose generic fibre has ample anticanonical divisor and the relative Picard rank $\rho(Y/B)$ is 1. It is rarely the case when there is a unique Mori fibre space representing its birational class and it is then a natural question to understand the relations between these. 

This survey takes this question in the case of weighted complete intersection Fano threefolds as its central theme. Already in this case, the geometry is so rich and mysterious that it has served as the theme of research for more than a few decades for many researchers. We hope that this survey will entice the community to continue doing so.

We work over the field of complex numbers.
\vspace{1cm}

\textit{Acknowledgments:} The authors would like to extend their heartfelt gratitute to Prof. Ivan Cheltsov for suggesting to us to write this survey. The first author was supported by ERC StG Saphidir No. 101076412.
The second author was supported by JSPS KAKENHI Grant Number 23K22389.

\section{Singularities}

\subsection{Terminal Singularities} \label{sect:termsing}

It turns out that in order to generalise the results on the birational classification of surfaces to higher dimensions, one needs to allow some mild singularities. See \cite[Section~16]{uenobook}.

\begin{Def}
We say that a normal projective variety $X$ is \textbf{$\mathbb{Q}$-factorial} if any Weil divisor $D$ on $X$ is $\mathbb{Q}$-Cartier, that is, there is a positive integer $m$ for which $mD$ is Cartier.
\end{Def}

In particular, we can take the pull-back of $D$ with respect to any morphism $W \rightarrow X$ and define the intersection number $D \cdot C$ with any curve $C \subset X$. Next we define the smallest class of singularities needed to run the Minimal Model Program.

\begin{Def}[Terminal Singularities]
Suppose $X$ is a normal projective variety. Then, $X$ has \textbf{terminal singularities} if 
\begin{itemize}
	\item $X$ is $\mathbb{Q}$-Gorenstein, that is, there is a positive integer $m$ for which the Weil divisor $mK_X$ is Cartier;
	\item if $f\colon Y \rightarrow X$ is a resolution of $X$ and $\{ E_i\}$ the collection of all exceptional prime divisors of $f$ then 
	\[
	K_Y = f^*(K_X) + \sum a_iE_i
	\]
	with strictly positive discrepancies $a_i >0$. 
\end{itemize}
We say that $X$ is in the \textbf{Mori Category} if it is a normal projective variety with $\mathbb{Q}$-factorial and terminal singularities. We say $X$ has \textbf{canonical singularities} if $a_i \geq 0$ for any $i$. 
\end{Def}

\begin{Rem} \label{rem:term}
    If a normal projective variety $X$ has terminal singularities then its singular locus has codimension at least $3$. See \cite[Corollary~5.18]{kollarmori}. In particular, terminal surfaces are smooth and terminal threefolds have isolated singularities.  
\end{Rem}

For higher dimensional varieties, the classification of terminal singularities is much more complicated and, indeed, widely open. In dimension three, however, the complete list of terminal singularities is known and is due to the works of Reid \cite{ReidI}, Danilov \cite{DanilovI}, Morrison and Stevens \cite{MorrisonStevens}, Mori \cite{MoriTerminal} and others. See \cite{reidyoung} for a friendly introduction to terminal and, more generally, canonical singularities.

\begin{Def}
Let $F \in \mathbb{C}\{x_1,x_2,x_3,x_4 \}$ be a convergent power series around $0$. Then $F=0$ is a \textbf{compound du Val singularity} (or cDV) if $F$ is of the form
\[
F = h(x,y,z)+tg(x,y,z,t)
\] 
where $h(x,y,z)=0$ defines a canonical surface singularity (also known as du Val singularity). 
\end{Def}

We have a more explicit description as follows. See \cite[Proposition~2.3.3]{kawakitabook} or \cite[Section~2]{kollarrealI}: 
\begin{Thm} \label{thm:cDV}
    Every cDV singularity is analytically given by one of the following convergent power series
\begin{table}[h]
  \centering

  \label{tab:arrow_table}
 \setlength{\tabcolsep}{10pt} 
  \renewcommand{\arraystretch}{1.2} 
  \begin{tabular}{@{} l l l  @{}} 
    \toprule
    \textbf{Type} & \textbf{ Function } &  \textbf{Conditions} \\
    \midrule
    $cA_{\leq 1}$ &  $x^2+y^2+z^2+h(t)$  & $h$ is not identically $0$ \\
     
    $cA_{\geq 2}$ & $x^2+y^2+h(z,t)$ &  $\mathrm{ord}\,h \geq 3$ \\

     $cD_{4}$ & $x^2+h(y,z,t)$  & $\mathrm{ord}\,h \geq 3$ and $h_3$ is not divisible by the square of a linear form. \\

      $cD_{\geq 5}$ & $ x^2+y^2z+ayt^r+h(z,t)$  & \makecell{\phantom{.} \\ \hspace{-5cm}$\mathrm{ord}\,h \geq s$, $a\in \mathbb C$, $r\geq 3,\, s\geq 4$ and $h_4\not =0$. \\ This is type $cD_n$ where $n=\min\{2r,s+1\}$ if $a\not =0$ and $n=s+1$ if $a=0$.}\\

      $cE_6$ & $x^2+y^3+yg(z,t)+h(z,t)$  & $\mathrm{ord}\,g \geq 3,\, \mathrm{ord}\,h \geq 4$ and $h_4\not =0$. \\

         $cE_7$ & $x^2+y^3+yg(z,t)+h(z,t)$  & $\mathrm{ord}\,g \geq 3,\, \mathrm{ord}\,h \geq 5$ and $g_3\not =0$. \\
         $cE_8$ & $x^2+y^3+yg(z,t)+h(z,t)$  & $\mathrm{ord}\,g \geq 4,\, \mathrm{ord}\,h \geq 5$ and $h_5\not =0$. \\
         
    \bottomrule
  \end{tabular}
    \caption{Terminal threefold singularities of index $1$.}
\end{table}
\end{Thm}
\begin{proof}
    You can see a proof in \cite[Theorem~2.8, Theorem~2.9, Theorem~2.10]{kollarrealI} or \cite[Proposition~2.3.3]{kawakitabook}.
\end{proof}



Recall that $\bm{\mu}_r$ denotes the cyclic group of $r$th roots of unity. Define the action of $\bm{\mu}_r$ on $\mathbb{C}^4$ with coordinates $x_1,\,x_2,\,x_3,\,x_4$ by
\begin{align*}
\bm{\mu}_r \times \mathbb{C}^4 &\longrightarrow \mathbb{C}^4 \\
(\epsilon, (x_1,x_2,x_3,x_4)) &\longmapsto (\epsilon^{\alpha_1}x_1,\epsilon^{\alpha_2}x_2,\epsilon^{\alpha_3}x_3,\epsilon^{\alpha_4}x_4)
\end{align*} 
where $\epsilon$ is a primitive $r$th root of unity and $\alpha_i$ are integers. We denote such an action by
\[
\frac{1}{r}(\alpha_1,\alpha_2,\alpha_3,\alpha_4)
\] 
and its orbit space by
\[
\mathbb{C}^4/\bm{\mu}_r(\alpha_1,\alpha_2,\alpha_3,\alpha_4).
\]
Let $F \in \mathbb{C}\{x_1,x_2,x_3,x_4 \}$ be a convergent power series around $0$ and suppose that $F$ is equivariant with respect to this action. Then $\bm{\mu}_r$ acts on the germ of the hypersurface $(F(x_1,x_2,x_3,x_4)=0) \subset \mathbb{C}^4$ and we can take the quotient
\[
(F(x_1,x_2,x_3,x_4)=0)/\bm{\mu}_r(\alpha_1,\alpha_2,\alpha_3,\alpha_4).
\]

The main theorem of \cite{reidterminal} is

\begin{Thm}[{\cite[Main Theorem~I]{reidterminal}}]
Every terminal 3-fold singularity over $\mathbb{C}$ is analytically isomorphic to
\[
(F(x_1,x_2,x_3,x_4)=0)/\bm{\mu}_r(\alpha_1,\alpha_2,\alpha_3,\alpha_4)
\]
where $F$ defines a cDV singularity.
\end{Thm}

It turns out that only very few actions produce terminal singularities. 
\begin{Thm}[{\cite{MoriTerminal}}]
    Let $0\in X$ be the germ of a terminal 3-fold analytic singularity over $\mathbb C$ of index $r>1$. Suppose that $0\in X$ is not a hypersurface singularity as in Theorem \ref{thm:cDV}. Then $(o\in X)$ is isomorphic to one of the following singularity types in Table \ref{tab:singindhigh}.
\begin{table}[h]
  \centering
 \setlength{\tabcolsep}{10pt} 
  \renewcommand{\arraystretch}{1.2} 
  \begin{tabular}{@{} l l l l @{}} 
    \toprule
    \textbf{Type} & \textbf{ Function } & \textbf{Action} & \textbf{Conditions} \\
    \midrule
    $cA/r$ &  $xy+g(z,t^r)$ & $\frac{1}{r}(1,-1,0,b)$ & $\mathrm{gcd}(b,r)=1$ \\
     
    $cAx/2$ & $x^2+y^2+g(z,t)$ & $\frac{1}{2}(1,0,1,1)$ & $\mathrm{ord}\,g \geq 4$ \\

     $cAx/4$ & $x^2+y^2+g(z,t^2)$ & $\frac{1}{4}(1,3,2,1)$ & $\mathrm{ord}\,g \geq 2$ \\

      $cD/2$ & $x^2+g(y,z,t)$ & $\frac{1}{2}(1,1,0,1)$ & $\mathrm{ord}\,g \geq 3,\, yzt\in g_3\,\, \text{or}\,\, y^2z\in g_3$\\

      $cD/3$ & $x^2+g(y,z,t)$ & $\frac{1}{3}(0,2,1,1)$ & $\mathrm{ord}\,g \geq 3,\, g_3\in\{y^3+z^3+t^3,y^3+zt^2,y^3+z^3\}$ \\

         $cE/2$ & $x^2+y^3+yh(z,t)+g(z,t)$ & $\frac{1}{2}(1,0,1,1)$ & $\mathrm{ord}\,h \geq 4,\,\,\mathrm{ord}\,g = 4$\\
    \bottomrule
  \end{tabular}
    \caption{Terminal threefold singularities of index greater than $1$.}
    \label{tab:singindhigh}
\end{table}

\end{Thm}

Of special importance for us is the case of \emph{cyclic quotient singularities}. These are actually $cA/r$ singularities where $g$ is linear. In this case, $F$ is the germ of a smooth point and 

$$
(o \in X) \simeq \frac{1}{r}(1,-1,b) \simeq  \frac{1}{r}(1,a,r-a)
$$
where $\mathrm{gcd}(b,r)=1$ or $\mathrm{gcd}(a,r)=1$. Indeed every threefold cyclic quotient singularity is terminal if and only if it is of the type above. See \cite[Section~5]{reidyoung}.

\subsection{Terminal Extractions} \label{sec:termextr} In the previous section we detailed the classification of terminal threefold singularities. In order to construct Sarkisov links it is essential to understand extremal divisorial contrations from various centres, including singular points.

\begin{Def}
Let $X$ be a normal variety (or a (analytic) germ of a normal variety) with only terminal singularities. In this survey, by a \textbf{divisorial
extraction} $\varphi \colon Y \rightarrow X$ centered along $\Gamma \subset X$, we mean a projective birational morphism from a normal variety with only terminal singularities such that $-K_Y$ is $\varphi$-ample, it contracts a prime divisor $E$ onto $\Gamma$ and it restricts to an isomorphism on $Y \setminus E$.
A divisorial extraction is called \textbf{extremal} if both $X$ and $Y$ are $\mathbb{Q}$-factorial projective varieties.
\end{Def}

Note that an extremal divisorial extraction $\varphi \colon Y \to X$ is nothing but a birational contraction of a $K_Y$-negative extremal ray from a variety $Y$ in the Mori category.
In order to construct Sarkisov links it is essential to understand extremal divisorial contractions, see \cite{cortiSP}.

The centre of a threefold divisorial extraction can be a point or a curve, possibly singular.

\subsubsection{The centre is a curve}  Terminal extractions from smooth curves in terminal $\mathbb Q$-factorial threefolds  were systematically analysed by Tziolas,  see \cite{tziolas2010three,tziolasI,tziolas2005three,tziolas2005families}. On the other hand, Ducat studied terminal extractions from singular curves in smooth threefolds, see \cite{ducatsingcurves}. See \cite[Chapter~4.2 and Chapter~4.3]{kawakitabook} for a discussion.

\subsubsection{The centre is a point} We collect known results when the centre is a point. Kawakita classified divisorial contractions to smooth points and showed that, local analytically, these are certain weighted blowups.

\begin{Thm}[{\cite[Theorem~2.2]{kawakitasmthpt}}]
    Let $\pi \colon E\subset Y \rightarrow P\in X$ be a threefold divisorial contraction which contracts the divisor $E$ to a smooth point $P\in X$. Then we can take local coordinates $x,\,y,\,z$ around $P$ and coprime positive integers such that $\pi$ is the weighted blowup of $P$ with weights $\mathrm{wt}(x,y,z)=(1,a,b)$.
\end{Thm}

The classification of terminal divisorial extractions from singular points has been carried out by Kawamata, Kawakita, Hayakawa and Yamamoto in several papers.

\begin{Thm}[{\cite{kawamata}}] \label{thm:kwbl}
Let $(\mathbf{p}\in X) \sim \frac{1}{r}(1,a,r-a)$ be the germ of a terminal cyclic quotient singularity. If $\varphi \colon Y \rightarrow X$ is a divisorial extraction centred along $\Gamma$ where $\mathbf{p} \in \Gamma$, then $\Gamma= \mathbf{p}$ and $\varphi$ is the weighted blowup with weights $\frac{1}{r}(1,a,r-a)$ and discrepancy $\frac{1}{r}$. Moreover, the exceptional divisor is $E \simeq \mathbb{P}(1,a,r-a) \subset Y$.
\end{Thm}
 
The weighted blowup of the previous theorem is called the \textbf{Kawamata blowup centred at $\mathbf{p}$}. An immediate corollary is the following

\begin{Cor}[{\cite[Corollary~3.4.3]{CPR}}] \label{cor:kwblcor}
Suppose $X$ is a 3-fold with only terminal quotient singularities. Suppose $\Gamma$ is a curve. If $\Gamma \subset X$ is the centre of a divisorial extraction in the Mori category then $\Gamma \subset X \setminus X_{\mathrm{sing}}$.
\end{Cor}

\begin{Thm}[{\cite[Theorem~1.1]{kawakitacA1}, \cite[Theorem~1.13]{kawakitaelephants}, \cite[Theorems~2.1 and 2.6]{yamamoto}, \cite[Theorem~6.1]{Erikcounting}}]
\label{thm:divcontcA}
Let $\mathbf{p} \in X$ be the germ of a $3$-fold terminal singularity of type $cA$ and let $\varphi \colon Y \to X$ be a divisorial extraction centred at $\mathbf{p} \in X$.
Then, one of the following holds.
\begin{enumerate}
\item There are positive integers $r_1, r_2, a$ and a power series $g (z, t) \in \mathbb{C} \{z, t\}$ such that the germ $\mathbf{p} \in X$ can be analytically identified with the germ
\[
o \in (x y + g (z, t) = 0) \subset \mathbb{C}^4,
\]
where $r_1 + r_2$ is divisible by $a$, $a$ is coprime to $r_i$ for $i = 1, 2$, $g (z, t)$ has weighted order $r_1 + r_2$ with respect to $\mathrm{wt} (z, t) = (a, 1)$ and $z^{(r_1+r_2)/a} \in g$.
Under the above identification, $\varphi$ is the weighted blowup with weights $\wt (x, y, z, t) = (r_1, r_2, a, 1)$ and discrepancy $a$.
\item The singularity $\mathbf{p} \in X$ is of type $cA_1$ and $\mathbf{p} \in X$ can be analytically identified with the germ
\[
o \in (x y + z^2 + t^3 = 0) \subset \mathbb{C}^4.
\]
After the above identification, $\varphi$ is the weighted blowup with weights $\wt (x, y, z, t) = (1, 5, 3, 2)$ and of discrepancy $4$.
\item The singularity $\mathbf{p} \in X$ is of type $cA_2$ and $\mathbf{p} \in X$ can be analytically identified with the germ
\[
o \in (x^2 + y^2 + z^3 + x t^2 = 0) \subset \mathbb{C}^4.
\]
After the above identification, $\varphi$ is the weighted blowup with weights $\wt (x, y, z, t) = (4, 3, 2, 1)$ and of discrepancy $3$.
\end{enumerate}
\end{Thm}

\begin{Thm}[{\cite[Theorem~1.1]{kawakitahigh}}]
Let $\mathbf{p} \in X$ be the germ of a $3$-fold terminal singularity of type $cA/r$ which is not a cyclic quotient singularity and let $\varphi \colon Y \to X$ be a divisorial extraction centred at $\mathbf{p} \in X$.
Then, the germ $\mathbf{p} \in X$ can be analytically identified with the germ
\[
o \in (x y + g (z^r, t) = 0)/\bm{\mu}_r (1, -1, 0, b),
\]
such that $\varphi$ is the weighted blowup with weights $\wt (x, y, z, t) = (r_1/r, r_2/r, a/r, 1)$ and of discrepancy $a/r$ which satisfies the following conditions.
\begin{enumerate}
\item $b$ is coprime to $r$.
\item $a \equiv b r_1 \pmod{r}$ and $r_1 + r_2 \equiv 0 \pmod{a r}$.
\item $(a-b r_1)/r$ is coprime to $r_1$.
\item $g (z^r, t)$ has weighted order $(r_1+r_2)/r$ with respect to the weights $\wt (z, t) = (a/r,1)$.
\end{enumerate}
\end{Thm}

\begin{Thm}[{\cite[Theorems~8.4 and 8.9]{HayakawaI}, \cite[Theorems~1.2 and 1.3]{kawakitahigh}}]
Let 
\[
\mathbf{p} \in X \cong o \in (x^2 + y^2 + g (z, t) = 0)/\bm{\mu}_2 (1, 0, 1, 1)
\]
be the germ of a $3$-fold terminal singularity of type $cAx/2$, where $g (z, t) \in (z, t)^4$ is $\bm{\mu}_2$-invariant.
Let $k$ be the weighted order of $g (z, t)$ with respect to the weights $\wt (z, t) = \frac{1}{2} (1, 1)$ and let $g_k (z, t)$ be the weight $d$ part of $g (z, t)$.
\begin{enumerate}
\item If $g_k (z, t)$ is not a square, then there is a unique divisorial contraction $\varphi \colon Y \to X$ centred at $\mathbf{p} \in X$, which is the weighted blowup with weights
\[
\wt (x, y, z, t) =
\begin{cases}
\frac{1}{2} (k+1, k, 1, 1), & \text{if $k$ is even}, \\
\frac{1}{2} (k, k+1, 1, 1), & \text{if $k$ is odd},
\end{cases}
\]
and discrepancy $1/2$.
\item If $g _k (z, t) = - h (z, t)^2$ for some $h (z, t) \in \mathbb{C} \{z, t\}$, then there are exactly $2$ divisorial extractions $\varphi_{\pm} \colon Y_{\pm} \to X$ centred at $\mathbf{p} \in X$.
The morphisms $\varphi_{\pm}$ is the composite of the weighted blowup $Y_{\pm} \to X_{\pm}$ of $\mathbf{p}_{\pm} \in X_{\pm}$, where
\[
\mathbf{p}_{\pm} \in X_{\pm} = 
\begin{cases}
o \in (x^2 + y^2 + \mp 2 y h (z, t) + g (z, t) - g_k (z, t) = 0)/\bm{\mu}_2 (1,0,1,1), & \text{if $k$ is even}, \\
o \in (x^2 \mp 2 x h (z, t) + y^2 + g (z, t) - g_k (z, t) = 0)/\bm{\mu}_2 (1,0,1,1), & \text{if $k$ is odd},
\end{cases}
\]
with weights
\[
\wt (x, y, z, t) =
\begin{cases}
\frac{1}{2} (k+1, k+2, 1, 1), & \text{if $k$ is even}, \\
\frac{1}{2} (k+2, k+1, 1, 1), & \text{if $k$ is odd},
\end{cases}
\]
and discrepancy $1/2$, and the isomorphism $\mathbf{p}_{\pm} \in X_{\pm} \cong \mathbf{p} \in X$ defined by $y \mapsto y \pm h (z, t)$ (or $x \mapsto x \pm h (z, t)$).
\end{enumerate}
\end{Thm}

\begin{Thm}[{\cite[Theorems~7.4 and 7.9]{HayakawaI}, \cite[Theorems~1.2 and 1.3]{kawakitahigh}}]
Let 
\[
\mathbf{p} \in X \cong o \in (x^2 + y^2 + g (z, t^2) = 0)/\bm{\mu}_4 (1, 3, 2, 1)
\]
be the germ of a $3$-fold terminal singularity of type $cAx/4$, where $g (z, t^2) \in (z, t)^2$ is $\bm{\mu}_4$-semi-invariant.
Let $(2k+1)/2$ be the weighted order of $g (z, t)$ with respect to the weights $\wt (z, t) = \frac{1}{4} (2, 1)$, where $k$ is a nonnegative integer, and let $g_{(2k+1)/2} (z, t)$ be the weight $(2k+1)/2$ part of $g (z, t)$.
\begin{enumerate}
\item If $g_{(2k+1)/2} (z, t)$ is not a square, then there is a unique divisorial contraction $\varphi \colon Y \to X$ centred at $\mathbf{p} \in X$, which is the weighted blowup with weights
\[
\wt (x, y, z, t) =
\begin{cases}
\frac{1}{4} (2k+1, 2k+3, 2, 1), & \text{if $k$ is even}, \\
\frac{1}{4} (2k+3, 2k+1, 2, 1), & \text{if $k$ is odd},
\end{cases}
\]
and discrepancy $1/4$.
\item If $g _{(2k+1)/2} (z, t) = - h (z, t)^2$ for some $h (z, t) \in \mathbb{C} \{z, t\}$, then there are exactly $2$ divisorial extractions $\varphi_{\pm} \colon Y_{\pm} \to X$ centred at $\mathbf{p} \in X$.
The morphisms $\varphi_{\pm}$ is the composite of the weighted blowup $Y_{\pm} \to X_{\pm}$ of $\mathbf{p}_{\pm} \in X_{\pm}$, where
\[
\mathbf{p}_{\pm} \in X_{\pm} = 
\begin{cases}
o \in (x^2 + \mp 2 x h (z, t) + y^2 + g (z, t) - g_{(2k+1)/2} (z, t) = 0)/\bm{\mu}_4 (1,3,2,1), & \text{if $k$ is even}, \\
o \in (x^2 + y^2 + \mp 2 y h (z, t) + g (z, t) - g_{(2k+1)/2} (z, t) = 0)/\bm{\mu}_4 (1,3,2,1), & \text{if $k$ is odd},
\end{cases}
\]
with weights
\[
\wt (x, y, z, t) =
\begin{cases}
\frac{1}{4} (2k+5, 2k+3, 2, 1), & \text{if $k$ is even}, \\
\frac{1}{4} (2k+3, 2k+5, 2, 1), & \text{if $k$ is odd},
\end{cases}
\]
and discrepancy $1/4$, and the isomorphism $\mathbf{p}_{\pm} \in X_{\pm} \cong \mathbf{p} \in X$ defined by $x \mapsto x \pm h (z, t)$ (or $y \mapsto y \pm h (z, t)$).
\end{enumerate}
\end{Thm}

The divisorial extractions centered at $3$-fold terminal singular points of types $cD/2$, $cD/3$ and $cE/2$ are also classified.
Because of the complexity of the classifications, we only give references to each singularity.
\begin{itemize}
\item The divisorial extractions with centre a $3$-fold terminal singular point of type $cD/2$ of discrepancy $1/3$ (resp.\ not equal $1/3$) are classified by \cite{HayakawaII} (resp.\ \cite[Theorems~1.2 and 1.3]{kawakitahigh}, \cite{kawakitaSuppl} and \cite{HayakawaII}).
\item It is proved in \cite[Theorems~1.2 and 1.3]{kawakitahigh} that any divisorial extraction with center a $3$-fold terminal singular point of type $cD/3$ is of discrepancy either $1/3$ and the divisorial extractions of discrepancy $1/3$ are classified in \cite[\S 9]{HayakawaI}.
\item It is proved in \cite[Theorems~1.2 and 1.3]{kawakitahigh} that any divisorial extraction with center a $3$-fold terminal singular point of type $cE/2$ is of discrepancy either $1/2$ or $1$.
Divisorial extractions of discrepancy $1/2$ and $1$ are classified in \cite[\S 10]{HayakawaI} and \cite[Theorem~1.2]{Hayakawadiscrep1}, respectively.
\end{itemize}

The remaining centers to consider are $3$-fold terminal singularities of types $cD$ and $cE$.
After the work of Kawakita \cite[Corollary~1.15]{kawakitaelephants}, Yamamoto \cite{yamamoto} completed the classification of divisorial extractions centered at a $3$-fold terminal singular point of types $cD$ and $cE$ of discrepancy greater than $1$.
The classifications of divisorial extractions centred at a $3$-fold singularity of types $cD$ and $cE$ of discrepancy $1$ are also completed in \cite{HayakawacD} and \cite{HayakawacE}, respectively, which are currently unpublished preprints and whose contents are knwon only to a very limited number of specialists.

\subsection{Small $\mathbb Q$-factorial modifications}

\begin{Def}
    A \textbf{small $\mathbb Q$-factorial modification}, or SQM, of a normal projective variety is a \emph{small} birational map $Y \dashrightarrow Y'$ to a normal $\mathbb Q$-factorial projective variety $Y'$, where a birational map $f \colon V \dashrightarrow V'$ between varieties is said to be \textbf{small} if it is an isomorphism in codimension one.
\end{Def}

\begin{Def}
    A \textbf{flipping contraction} $\pi \colon Y \rightarrow Z$ is a small extremal contraction from a terminal $\mathbb Q$-factorial variety such that $-K_Y$ is $\pi$-ample. The \textbf{flip} of $\pi$ is a small extremal contraction $\pi^+\colon Y^+\rightarrow Z$ such that $K_{Y^+}$ is $\pi^+$-ample. We usually represent a flip as a diagram
      \[
        \begin{tikzcd}[ampersand replacement=\&,column sep = 1.5em]
             \&   Y  \ar[rr,dashed,"\tau"] \ar[rd,swap,"\pi"]  \& \& Y^+  \ar[ld,"\pi^+"]  \&  \\
               \& \& Z \&  \& 
        \end{tikzcd}
    \]
    and call $\tau \colon Y \dashrightarrow Y^+$ the flip of $Y$ by abuse of notation. A \textbf{flopping contraction} and a \textbf{flop} can be defined similarly except that we require $K_{Y}$ (resp. $K_{Y^+}$) to be numberically $\pi$-trivial
 (resp.  numerically $\pi^+$-trivial). If $\tau$ is the flip of $Y$, an \textbf{antiflip} is the inverse map $Y^+\dashrightarrow Y$.
 \end{Def}

 \begin{Rem}
 Mori proved in the monumental paper \cite{mori1988flip} the existence of threefold flips. See also \cite{KMflips}. The existence of flips was settled in higher dimensions in \cite{BCHM}.
 \end{Rem}

\begin{Rem}
    If $\tau \colon Y \dashrightarrow Y^+$ is a flip, then $Y^+$ is normal $\mathbb Q$-factorial and $\rho(Y)=\rho(Y^+)$. See \cite[Lemma~1.3.20]{kawakitabook}. In particular, a flip is an example of an SQM.
\end{Rem}

\begin{Rem}
    The flip or flop of a terminal variety is again terminal. See \cite[Corollary~3.42]{kollarmori}. In the case of three dimensional flops the analytic type of the singularities is actually preserved. See \cite[Theorem~6.15]{kollarmori}. The assumption on the dimension is essential. See, for instance, \cite[Example~3.18]{campo2023blowups} for an example of a flop between a smooth fivefold and a singular one. On the other hand, suppose that a variety $Y$ in the Mori category admits a small contraction $\pi \colon Y \to Z$ of a $K_Y$-positive extremal ray and suppose further that there is a small $K_{Y^+}$-negative extremal contraction $\pi^+ \colon Y^+ \to Z$. In this case the variety $Y'$ may not be terminal.
\end{Rem}

\begin{Ex}
    Consider the vector bundle 
    $\mathcal{E}=\mathcal{O}_{\mathbb P^1}(-2)\oplus\mathcal{O}_{\mathbb P^1}\oplus\mathcal{O}_{\mathbb P^1}$ and $T=\mathbb P(\mathcal E)$, which is  a smooth toric variety. Besides the identity $T$ has an SQM to $T'$ which is a fibration over $\mathbb P^1$ and the fibre over each point is a cone over a conic, so that $T'$ has a line of singularities and therefore is not terminal. See Remark \ref{rem:term}.
\end{Ex}

\begin{Ex}
    On the other hand, there are terminal antiflips. Consider the vector bundle 
    $\mathcal{E}=\mathcal{O}_{\mathbb P^1}(-2)\oplus\mathcal{O}_{\mathbb P^1}(-1)\oplus\mathcal{O}_{\mathbb P^1}$ and $T=\mathbb P(\mathcal E)$, which is  a smooth toric variety. Besides the identity $T$ the only other SQM is an antiflip to a toric variety $T'$ which has a unique $\frac{1}{2}(1,1,1)$ singularity. Hence $T'$ is terminal.
\end{Ex}

\section {Weighted Complete Intersections}
\begin{Def}
\label{Def:Fano}
A normal projective variety over $\mathbb{C}$ is called a \textbf{Fano $n$-fold} if $X$ is an $n$-dimensional variety with $\mathbb{Q}$-factorial terminal singularities and such that $-K_X$ is ample. The \textbf{Fano index} of $X$ is 
\[
\iota_X := \max \{q \in \mathbb{Z}_{>0} \, | \, -K_X \sim_{\mathbb{Q}} qA,\, \text{$A \in \Cl(X)$ is a $\mathbb{Q}$-Cartier Divisor}\}.
\]  
If $\iota_X \geq 2$ we say that $-K_X$ is \textbf{divisible} in $\Cl(X)$ and that $X$ has \textbf{higher index}.
\end{Def}

Let $A \in \Cl(X)$ be an ample divisor for which $-K_X \sim_{\mathbb{Q}} \iota_XA$. We consider $X$ to be polarised by $A$, that is, we consider $X$ together with an embedding into weighted projective space given by the ring of sections,
\[
R(X,A):=\bigoplus_{m \geq 0} H^0(X,\mathcal{O}_X(mA)).
\] 
Notice that $R(X,-K_X) \subset R(X,A)$. These are finitely generated over $\mathbb{C}$, see \cite[Proposition~1.11(4)]{mdsGIT}, and we realise $X$ as
\[
X:=\Proj R(X,A).
\]
A choice of generators $x_0,\ldots,x_N$ where $x_i \in H^0(X,\mathcal{O}_X(a_iA))$ determines a map 
\[
X \hookrightarrow \Proj\mathbb{C}[x_0, \ldots,x_N] = \mathbb{P}(a_0,\ldots,a_N)=\mathbb{P}.
\]  For a choice of a minimal number of generators $x_0, \dots, x_N$, we define the \textbf{codimension} of the 3-fold $X$, denoted by $\mathrm{cod}_X$, to be $N-3$. A subvariety $X \subset \mathbb{P}$ of codimension $N-3$ is a \textbf{weighted complete intersection (WCI) of multidegree $(d_1,\ldots, d_{N-3})$} if its weighted homogeneous ideal $I_X \subset \mathbb{C}[x_0,\ldots,x_N]$ is generated by a regular sequence of homogeneous elements of degrees $d_i$ in  $\mathbb{C}[x_0,\ldots,x_N]$.  Unless otherwise stated, we always assume that $\mathbb{P}(a_0,\ldots,a_N)$ is \textbf{well formed}, that is, we have $\gcd(a_0,\ldots, \widehat{a_i},\ldots, a_N) = 1$ for each $i$. Each weighted projective space is isomorphic to a well formed weighted projective space, see \cite[1.3.1]{DolgachevWPS}. A subvariety $X \subset \mathbb{P}$ is \textbf{well formed} if $\mathbb{P}$ is well formed and $\codim_X(X \cap \mathbb{P}_{\text{sing}}) \geq 2$.

\begin{Def}[{\cite[Definition~3.1.5]{DolgachevWPS}}]
Consider the weighted projective space $\mathbb{P}$ and the natural quotient map $\pi \colon \mathbb{A}^{n+1} \setminus 0 \rightarrow \mathbb{P}$. A subvariety $X \subset \mathbb{P}$ is said to be \textbf{quasismooth} if its affine cone $C_X:=\pi^{-1}X$ is smooth outside the origin.
\end{Def}

\begin{Thm}[{\cite[Theorem~3.2.4]{DolgachevWPS}}] \label{thm:pic1}
Assume $X \subset \mathbb{P}$ is a quasismooth and well formed WCI of dimension at least $3$. Then $\Pic(X)=\mathbb{Z}$.
\end{Thm}

When $X$ is quasismooth its singularities are exclusively due to the $\mathbb{C}^{*}$-action defining the weighted projective space $\mathbb{P}$. If $X$ is additionally terminal, its singularities are terminal cyclic quotient singularities.

Throughout the survey, a Fano threefold weighted hypersurface (or more generally, WCI) is always $\mathbb{Q}$-factorial and has only terminal singularities even if we do not mention its singularities.
In particular, a quasismooth Fano threefold WCI has only terminal quotient singularities.
This is simply because a Fano threefold WCI is a Fano variety and a Fano variety is by definition $\mathbb{Q}$-factorial and has terminal singularities.

\begin{Rem}
Kawamata proves in \cite[Theorem~2]{Kawamatabound} that terminal Fano threefolds form a bounded family.   
\end{Rem}

\section{The various forms of rigidity}

\subsection{Definitions of rigidity}

Let $W$ be a smooth projective variety. If $K_W$ is not pseudo-effective, then by \cite[Corollary~1.3.3]{BCHM}, the Minimal Model Program produces a birational model $V$ of $W$ which admits a very special fibre structure $V \rightarrow S$ called a Mori fibre space.

\begin{Def}
A \textbf{Mori fibre space} is an extremal contraction $f \colon V \rightarrow S$ between normal projective varieties for which $f_*\mathcal{O}_V=\mathcal{O}_S$ and
\begin{itemize}
	\item $V$ has at worst $\mathbb{Q}$-factorial, terminal singularities;
	\item $-K_V$ is relatively ample for $f$;
	\item The relative Picard rank satisfies $\rho (V/S)=1$ and $\dim S < \dim V$.
\end{itemize}
We sometimes denote the Mori fibre space $f \colon V \rightarrow S$ by $V/S$ or simply by $V$ if the base is understood from the context. If $X$ is a Fano $n$-fold with $\rho(X)=1$, then we can consider it as a Mori fibre space where the base is a point. A Mori fibre space is \textbf{strict} if $\mathrm{dim} S>0$, in which case the Picard rank of the total space, $V$, is greater than $1$. 
\end{Def}

\begin{Ex}
Suppose $X$ is a quasismooth well-formed WCI of dimension $n\geq 3$. Suppose additionally that the singularities of $X$ are terminal. Then it follows from Theorem \ref{thm:pic1} that $X$ is a Mori fibre space.    
\end{Ex}

\begin{Ex}
    The Hirzebruch surfaces $F_n=\mathbb P(\mathcal{O}_{\mathbb P^1}\oplus\mathcal {O}_{\mathbb P^1}(-n)) \rightarrow \mathbb P^1$ are examples of two dimensional Mori fibre spaces.
\end{Ex}

\begin{Rem}
    The base of a Mori fibre space is also $\mathbb Q$-factorial (See \cite[Lemma~1.4.13]{kawakitabook}), but it is not necessarily terminal. However, the singularities of the base are not arbitrary. Indeed, when $V$ is a threefold and $S$ is a surface, the base can have at most du Val singularities of type $A$. See \cite[Theorem~1.2.7]{MoriProkConicsI}.  
\end{Rem}

\begin{figure*}[h!]
        \centering
        \begin{subfigure}[b]{0.475\textwidth}
            \centering
\begin{tikzcd}
Y\ar[d] \ar[rr, dashed]&  &  X' \ar[d, "f'"]    \\
X \ar[d, swap,"f"]  & & S'\ar[dll]\\
S &  &    
\end{tikzcd}

            \caption[Link1]%
            {{\small  Type I}}    
            \label{fig:1}
        \end{subfigure}
        \hfill
        \begin{subfigure}[b]{0.475\textwidth}  
            \centering 
\begin{tikzcd}
Y\ar[d] \ar[rr, dashed]&  &  Y' \ar[d]    \\
X \ar[d, swap,"f"]  & & X'\ar[d,"f'"]\\
S \ar[rr,equal] &  & S'    
\end{tikzcd}
            \caption[]%
            {{\small Type II}}    
            \label{fig:2}
        \end{subfigure}
        \vskip\baselineskip
        \begin{subfigure}[b]{0.475\textwidth}   
            \centering 
\begin{tikzcd}
X\ar[d, "f",swap] \ar[rr, dashed]&  &  Y' \ar[d]    \\
S \ar[drr]   & & X'\ar[d, "f'"]\\
 &  & S'   
\end{tikzcd}
            \caption[]%
            {{\small Type III}}    
            \label{fig:3}
        \end{subfigure}
        \hfill
        \begin{subfigure}[b]{0.475\textwidth}   
            \centering 
\begin{tikzcd}
X\ar[d,swap, "f"] \ar[rr, dashed]&  &  X' \ar[d, "f'", dashed]    \\
S \ar[dr]  & & S'\ar[dl]\\
 & R &    
\end{tikzcd}
            \caption[]%
            {{\small Type IV}}    
            \label{fig:4}
        \end{subfigure}
        \caption{\small The four types of Sarkisov links between Mori fibre spaces.} 
        \label{fig:4types}
    \end{figure*}
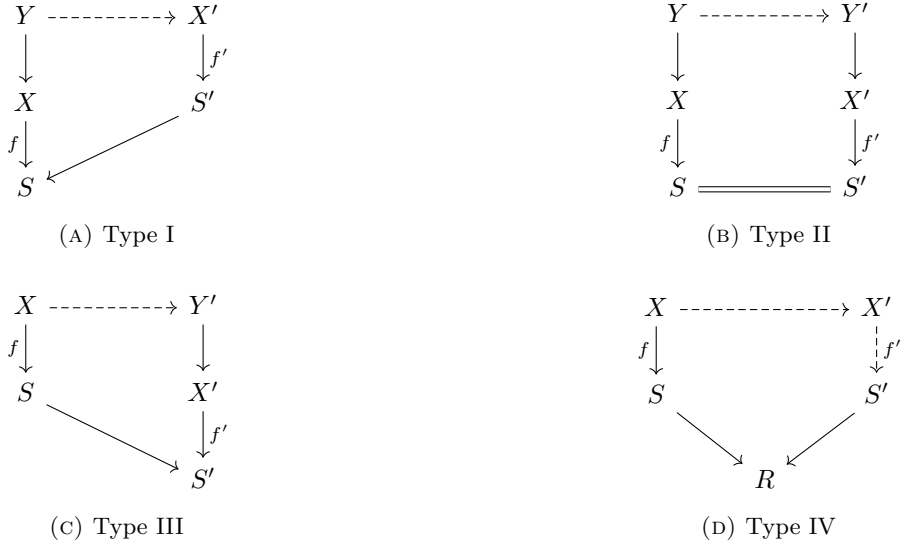

 The Sarkisov Program is the study of the relations between birational Mori fibre spaces. The following is the fundamental result of the Sarkisov Program: 

\begin{Thm}[{\cite[Theorem~1.1]{genSP}, \cite{cortiSP} for $3$-folds}] \label{thm:sp}
Two Mori fibre spaces are birational if and only if they are connected by a finite sequence of Sarkisov links.
\end{Thm}
A \textbf{Sarkisov link} $\sigma \colon X \dashrightarrow X'$ between two Mori fibre spaces $f \colon X \rightarrow S$ and $f' \colon X' \rightarrow S'$ is one of four types as in figure \ref{fig:4types}, see \cite{cortiSP}.

Every non-horizontal arrow is an extremal contraction. The dashed arrows are (finite) sequences of anti-flips, a flop or flips, in that order. Notice that a Sarkisov link from a Fano variety can only be of Type I or Type II, and it is always initiated by a divisorial contraction $Y \to X$. 

\begin{Ex}
    We can exhibit examples of Sarkisov links of all four types already in the case of surfaces. Let $F_n=\mathbb P(\mathcal{O}_{\mathbb P^1}\oplus\mathcal {O}_{\mathbb P^1}(-n)) \rightarrow \mathbb P^1$ be the $n$-th Hirzebruch surface. The blowup of $\mathbb P^2$ at a point is the $1$-st Hirzebruch Surface. The induced map $\mathbb P^2 \dashrightarrow F_1$ is a Sarkisov link of type I and the inverse is a Sarkisov link of type III. The surface $F_0$ is isomorphic to $\mathbb P^1 \times \mathbb P^1$ and it has two projections to $\mathbb P^1$. The isomorphism between the two projective bundles is a Sarkisov link of type IV with $R=\mathrm{Spec}\, \mathbb C$.

    Finally let $\widetilde{F_n}$ be the blowup of $F_n$ at the point $p$ in the section $\sigma_n = \mathbb P (\mathcal O_{\mathbb  P^1} (-n))$ of $F_n$. Then the fibre passing through $p$ becomes a $(-1)$-curve $\tilde f$ in $\widetilde{F_n}$ and the strict transform of $\sigma_n$ is a curve of self-intersection $-(n+1)$. By Castelnuovo's Contraction theorem, $\tilde f$ is the exceptional divisor of the blowup of a smooth surface $S$. But $S$ is a rational ruled surface with a unique section of self-intersection $-(n+1)$. Hence $S$ is $F_{n+1}$. The induced map $F_n \dashrightarrow F_{n+1}$ is a Sarkisov link of type II.
 \end{Ex}

The following two definitions capture the extent up to which a Fano variety of Picard rank $1$ is unique among Mori fibre spaces. The first one is essentially due to Corti (See \cite[Definition~1.3]{cortising} and \cite[Definition~1.1.2]{chel}) and the next one due to Shokuhrov and later by Abban-Okada, (See \cite[Definition~1.4]{hamidokadapfaff}).
\begin{Def}
A Fano variety $X$ of Picard rank $1$ is \textbf{birationally rigid} if whenever there is a birational map $\sigma \colon X \rat X'$ to a Mori fibre space $X'/S$, then $X$ and $X'$ are biregularlly isomorphic (and in particular $X'/S$ is a Fano variety of Picard rank $1$). If, in addition, the birational automorphism group of $X$ coincides with its biregular automorphism group we say that $X$ is \textbf{birationally superrigid}.  
\end{Def}

\begin{Def}[{\cite[Definition~1.4]{hamidokadapfaff}}]
    A Fano variety of Picard rank $1$ is \textbf{birationally solid} if it is not birational to any strict Mori fibre space.
\end{Def}

\begin{Rem}
    Clearly we have the following implications
    $$
    \text{Birationally superrigid} \quad \implies \quad \text{Birationally rigid} \quad \implies \quad \text{Birationally solid}  \quad \implies\quad \text{Irrational}
    $$
    However, none of the reverse implications hold.
    A factorial quartic threefold with only nodes (i.e.\ ordinary double points) is birationally rigid \cite[Theorem~2]{Mella} and the blowup at a node initiates a birational self-link \cite[\S~1]{Mella}.
    This provides an example of birationally rigid but not superrigid Fano varieties.
    A general quartic $3$-fold admitting a unique singular point (of type $cA_2$) which is analytically equivalent to $o \in (xy + z^3 + t^3 = 0) \subset \mathbb{C}^4$ is birationally solid but not birationally rigid (see \cite{cortimella}).
    It is well known that a smooth cubic threefold $X$ is irrational \cite{CG}, whereas the blowup $Y \to X$ along a line on $X$ resolves the indeterminacy of the projection $X \dashrightarrow \mathbb{P}^2$ and obtain a conic bundle $Y \to \mathbb{P}^2$.
    This provides an example of irrational but not birationally solid Fano varieties. See also \cite[Main Theorem 2]{Sarikyan} for yet another interesting example of a Fano threefold which has $8$ birational models of del Pezzo fibrations of degree $1$ over $\mathbb{P}^1$ and it is not birational to any other Mori fibre space (hence it is irrational).
\end{Rem}

\subsection{Behaviour in families}
For some time it was expected that birational rigidity was an open property in moduli, see \cite[Conjecture~1.4]{cortising}. This expectation was shown to be false in \cite[Corollary1.15]{rignotopen} by showing the existence of a birationally rigid complete intersection of a quadric and a cubic admitting a small deformation to a non-birationally rigid Fano threefold. Birational superrigidity is also not an open property in moduli, see \cite[Example~6.3]{chel}. On the other hand, it is still unknown whether birational solidity in an open property in moduli. We leave here the precise expectation:

\begin{Quest}[{\cite[Question~5.2]{zhuanglocallyclosed}}]
    Let $f \colon \mathcal X \to T$ be a $\mathbb Q$-Gorenstein family of Fano varieties. Is the set
    $$
    \{t\in T \,\, |\, \, \mathcal X_t \, \text{is birationally solid}\}
    $$
    an open subset of $T$ ?
\end{Quest}

\subsection{Maximal singularities}

The existence of a non-biregular birational map $X \dashrightarrow X'$ implies the existence of a mobile linear system on $X$ with a nonempty base locus.
We can capture the failure of birational (super)rigidity by the existence of a mobile linear system that is highly singular.
The key observation here is that any Sarkisov link from a Fano variety of Picard rank $1$ is initiated by a divisorial extraction.

Let $X$ be a normal variety, $\mathcal{M}$ a linear system of $\mathbb{Q}$-Cartier Weil divisors on $X$ and $E$ a prime divisor over $X$, that is, $E$ is a prime divisor on a normal variety $Y$ admitting a projective birational morphism $\varphi \colon Y \to X$.
We define $\ord_E (\mathcal{M})$ to be the number $\ord_E (\varphi^*D)$, where $D \in \mathcal{M}$ is a general member.
We say that a pair $(X, \lambda \mathcal{M})$, where $\lambda > 0$ is a rational number, is \textbf{canonical at $E$} if the inequality
\[
a_X (E) \ge \lambda \ord_E (\mathcal{M}),
\]
holds, where $a_X (E) := \ord_E (K_Y - \varphi^*K_X)$ is the discrepancy of $K_X$ at $E$.
We say that $(X, \lambda \mathcal{M})$ is \textbf{canonical} if it is canonical at any prime divisor over $X$.

\begin{Def}
Let $X$ be a Fano variety with $\rho (X) = 1$.
A divisorial extraction $\varphi \colon Y \to X$ is said to be a \textbf{maximal extraction} if there is a mobile linear system $\mathcal{M}$ on $X$ such that the pair $(X, \frac{1}{n} \mathcal{M})$, where $n > 0$ is the rational number determined by $\mathcal{M} \sim_{\mathbb{Q}} - n K_X$, is not canonical at the exceptional divisor of $\varphi$.
The centre of a maximal extraction is called a \textbf{maximal centre}.
\end{Def}

The following gives a characterization of birational superrigidity in terms of maximal centers/extractions. 

\begin{Thm}[{Noether--Fano--Iskovskikh inequality, \cite{iskmanin}, \cite[Proposition-definition 2.10]{cortiSP}, \cite[Theorem~1.26]{chelshra}}]
\label{thm:NFI}
Let $X$ be a Fano variety with $\rho (X) = 1$.
Then $X$ is birationally superrigid if and only if there is no maximal center on $X$.
\end{Thm}

\begin{proof}
If $X$ is not birationally superrigid, then by \cite[Theorem~1.26]{chelshra} (or by the argument of \cite{iskmanin}) there is a mobile linear system $\mathcal{M} \sim_{\mathbb{Q}} - n K_X$ such that the pair $(X, \frac{1}{n} \mathcal{M})$ is not canonical.
It then follows from \cite[Proposition-definition 2.10]{cortiSP} that there is a maximal extraction.
The converse follows from \cite[Theorem~1.26]{chelshra}.
\end{proof}

The difference between birational superrigidity and rigidity is the existence of a non-biregular birational automorphism of $X$. 

\begin{Thm}
Let $X$ be a Fano variety with $\rho (X) = 1$.
Then the following are equivalent.
\begin{enumerate}
\item $X$ is birationally rigid.
\item For each divisorial extraction $\varphi \colon Y \to X$, either it is not a maximal extraction or there is a Sarkisov self-link initiated by $\varphi$.
\end{enumerate}
\end{Thm}

\begin{proof}
Suppose $X$ is birationally rigid.
If it is superrigid, then (2) holds in a trivial way by Theorem~\ref{thm:NFI}.
If not, then there is a birational automorphism of $X$ which is not biregular.
By the Sarkisov program, such a birational automorphism is the composite of Sarkisov links.
Any link appearing in the decomposition must be a Sarkisov self-link since $X$ is birationally rigid.
In particular, there is a Sarkisov self-link of $X$ and the divisorial extraction initiating the self-link is a maximal extraction by \cite[Lemma~2.5]{okadathreeBMF}.

Suppose that $X$ is not birationally rigid.
Then there is a Sarkisov link $X \dashrightarrow X'$ to a Mori fibre space $X'/S'$ which is not isomorphic to $X$.
Again by \cite[Lemma~2.5]{okadathreeBMF}, the divisorial extraction initiating the link is a maximal center, and thus (2) does not hold.
\end{proof}

It follows that the proof of birational rigidity of a given Fano variety can be divided into two parts: (i) exclusion of maximal centers/extractions and (ii) construction of Sarkisov self-links.
In \S~\ref{sec:exclusion} and \S~\ref{sec:birinv}, we will explain (i) and (ii), respectively.

\subsection{Various methods of exclusion}
\label{sec:exclusion}

There are various methods for excluding a subvariety or a divisorial extraction as a maximal center or a maximal extraction.
We explain some of them.

A typical argument for excluding smooth points is as follows: Suppose a smooth point on a Fano threefold is a maximal center of a mobile linear system. Then the linear system has a base locus at the point with high multiplicity. The exclusion will be complete if this leads to a contradiction.
This method is based on the local inequality for intersection cycles. 

\begin{Thm}[{$4n^2$-inequality, \cite[Theorem~5.3.2]{CPR}}]
\label{thm:4nineq}
Let $\mathbf{p} \in X$ be a germ of a smooth threefold, $\mathcal{M}$ a mobile linear system on $X$ and $n > 0$ a rational number.
If $\mathbf{p}$ is the center of non-canonical singularity of the pair $(X, \frac{1}{n} \mathcal{M})$, then for general members $D_1, D_2 \in \mathcal{M}$, we have
\[
\mult_{\mathbf{p}} (D_1 \cdot D_2) > 4n^2.
\]
\end{Thm}

In the following, let $X$ be a Fano threefold with $\Cl (X) \cong \mathbb{Z}$ and let $A \in \Cl (X)$ be the positive generator of $\Cl (X)$. 

\begin{Lem}[{\cite[Lemma~2.16, Corollary~2.17]{okadaII}}]
\label{lem:exclNE}
Let $\varphi \colon Y \to X$ be a divisorial extraction.
If $\varphi$ is a maximal extraction, then the $1$-cycle $(-K_Y)^2$ is in the interior of the cone $\overline{\mathrm{NE}} (Y)$ of effective curves on $Y$.
In particular, if there exists a nef divisor $N$ on $Y$ such that $N \cdot (-K_Y)^2 \leq 0$, then $\varphi$ is not a maximal extraction.
\end{Lem}

The excluding method by fining a nef divisor $N$ as in Lemma~\ref{lem:exclNE} is called a test class method and such a nef divisor is called a test class.
We introduce the notion of isolating class which is useful when we find a test class and also when we derive a contradiction by applying Theorem \ref{thm:4nineq}.

\begin{Def}
Let $L$ be a Weil divisor class on $X$ and $\Gamma \subset X$ an irreducible subvariety.
For a positive integer $m$, we consider a linear system $\mathcal{L}_{\Gamma}^m := |\mathcal{I}_{\Gamma}^m (mL)|$, where $\mathcal{I}_{\Gamma}$ is the ideal sheaf of $\Gamma \subset X$. 
We say that the class $L$ \textbf{isolates} $\Gamma$ (or $L$ is a $\Gamma$-\textbf{isolating class}) if there is a positive integer $m$ such that $\Gamma$ is an isolated component of the linear system $\mathcal{L}_{\Gamma}^m$ and, in the case $\Gamma$ is a curve, the generic point of $\Gamma$ appears in $\Bs \mathcal{L}_{\Gamma}^m$ with multiplicity $1$.
\end{Def}

The following are methods for excluding curves.

\begin{Lem}[{\cite[Lemma 2.9]{okadaII}}]
\label{lem:exclCurvebydeg}
Let $C \subset X$ be an irreducible and reduced curve.
If $-K_X \cdot C \geq -K_X^3$, then $C$ is not a maximal centre. 
\end{Lem}

\begin{Lem}[{\cite[Lemma~2.11]{okadaII}}]
\label{lem:exclCurveint}
Let $\Gamma \subset X$ be an irreducible and reduced curve.
Assume that there is an effective divisor $S$ on $X$ containing $\Gamma$ and a movable linear system $\mathcal{M}$ on $X$ whose base locus contains $\Gamma$ satisfying the following properties.
\begin{enumerate}
\item $S \sim_{\mathbb{Q}} - m K_X$ for some rational number $m \ge 1$.
\item For a general member $T \in \mathcal{M}$, $T$ is a normal surface, the intersection $S \cap T$ is contained in the base locus of $\mathcal{M}$ set-theoretically, and the scheme $S \cap T$ is reduced along $\Gamma$.
\item Let $T \in \mathcal{M}$ be a general member and let $\Gamma, \Gamma_1, \dots, \Gamma_l$ be the irreducible and reduced curves contained in the base locus of $\mathcal{M}$.
For each $i = 1, \dots, l$, there is an effective $1$-cycle $\Delta_i$ on $T$ such that $\Gamma \cdot \Delta_i \geq -K_X \cdot \Delta_i > 0$ and $\Gamma_j \cdot \Delta_i \geq 0$ for $j \ne i$.
\end{enumerate}
Then, $\Gamma$ is not a maximal centre.
\end{Lem}

When we attempt to exclude a given divisorial extraction as a maximal extraction, the test class method or its variants is a first choice to apply.
If this fails, then there are some other methods. The following is one of such methods.

\begin{Lem}[{\cite[Lemma~2.20]{okadaII}}]
\label{lem:exclinfiniteC}
Let $\varphi \colon Y \to X$ be a divisorial extraction and let $E$ be its exceptional divisor.
Suppose that there are infinitely many curves $C_{\lambda}$ on $Y$ such that $-K_Y \cdot C_{\lambda} \leq 0$ and $E \cdot C_{\lambda} > 0$.
Then $\varphi$ is not a maximal extraction.
\end{Lem}

\subsection{Birational involutions}
\label{sec:birinv}

We explain some typical constructions of Sarkisov self-links of Fano varieties of Picard rank $1$.

Quadratic involutions are Sarkisov self-links of Fano varieties that can be constructed in the following way.

\begin{Prop}[Quadratic involutions, {\cite[Lemma 3.2]{okadaII}}] \label{prop:QI}
Let $X$ be a Fano variety with $\rho (X) = 1$ and let $\varphi \colon Y \to X$ be a divisorial extraction.
Assume that $-K_Y$ is nef and big but not ample, and there is a double cover $\pi \colon Z \to V$ from the anticanonical model $Z$ of $Y$ to a normal projective $\mathbb{Q}$-factorial variety $V$ with $\rho (V) = 1$.
Let $\iota_Z \colon Z \to Z$ be the biregular involution which swaps the fibres of $\pi$, and let $\iota_Y \colon Y \dashrightarrow Y$ be the induced birational involution, where $\psi \colon Y \to Z$ is the anticanonical morphism.
Then one of the following holds.
\begin{enumerate}
    \item $\psi$ is divisorial. In this case $\varphi$ is not a maximal extraction.
    \item $\psi$ is small. In this case $\iota_Y$ is the flop of $\psi$ and the diagram
    \begin{equation} \label{eq:QI}
    \vcenter{
    \xymatrix{
    Y \ar[d]_{\varphi} \ar[rd]^{\psi} \ar@{-->}[rrr]^{\iota_Y} & & & \ar[ld]_{\psi} Y \ar[d]^{\varphi} \\
    X & Z \ar[r]_{\iota_Z} & Z & X
    }
    }
    \end{equation}
    gives a Sarkisov self-link $\iota := \varphi \circ \iota_Y \circ \varphi^{-1} X \dashrightarrow X$.
\end{enumerate}
\end{Prop}

\begin{Thm}[{\cite[\S 4.1]{chel}, \cite[\S 4.4]{CPR}}]
Let $X = X_d \subset \mathbb{P} (a_{i_1}, a_{i_2}, a_{i_3}, a_j, a_k)$ be a quasismooth Fano threefold weighted hypersurface of index $1$ and of degree $d$, where $\{i_1, i_2, i_3, j, k\} = \{0, 1, 2, 3, 4\}$.
Assume that $x_k^2 x_j$ appears in the defining polynomial of $X$ with nonzero coefficient.
Then, after replacing homogeneous coordinates, the defining equation of $X$ can be written as
\[
x_k^2 x_j + x_k f + g = 0,
\]
where $f, g \in \mathbb{C} [x_{i_1}, x_{i_2}, x_{i_3}, j]$ are homogeneous polynomials of degree $a_k + a_j, 2 a_k + a_j$, respectively.
Let $\mathbf{p} \in X$ be the point at which only the homogeneous coordinate $x_k$ does not vanish and let $\varphi \colon Y \to X$ be the Kawamata blowup at $\mathbf{p}$.
\begin{enumerate}
\item If $f$ is divisible by $x_j$, then $\varphi$ is not a maximal extraction.
\item If $f$ is not divisible by $x_j$, then $\varphi$ initiates a Sarkisov self link $\iota \colon X \dashrightarrow X$ which is a birational involution.
\end{enumerate}
\end{Thm}

\begin{proof}
By the table of the $95$ families, we have $a_{i_1}, a_{i_2}, a_{i_3} < a_k$.
The projection $X \dashrightarrow \mathbb{P} := \mathbb{P} (a_{i_1}, a_{i_2}, a_{i_3}, a_j)$ to the coordinates $x_{i_1},x_{i_2},x_{i_3},x_j$ is a generically finite map of degree $2$.
The Kawamata blowup $\varphi \colon Y \to X$ at $\mathbf{p}$ is the weighted blowup with weights $\wt (x_{i_1}, x_{i_2}, x_{i_3}) = \frac{1}{a_k} (a_{i_1}, a_{i_2}, a_{i_3})$ since $a_{i_1}, a_{i_2}, a_{i_3} < a_k$, and it resolves the indeterminacy of $X \dashrightarrow \mathbb{P}$. 
We denote by $Y \xrightarrow{\psi} Z \xrightarrow{\pi} \mathbb{P}$ be the Stein factorization of the induced morphism $Y \to \mathbb{P}$, where $\psi$ is birational and $\pi$ is finite of degree $2$.

We claim that $\psi$ is the anticanonical morphism.
To see this, it is eanough to show that $x_{i_1}, x_{i_2}, x_{i_3}$ and $x_j$ lift to plurianticanonical sections on $Y$.
This is clear for $x_{i_1}, x_{i_2}, x_{i_3}$ since $x_{i_l}$ lift to a section of $-a_{i_l} K_Y = a_{i_l} \varphi^*A - \frac{a_{i_l}}{a_k} E$ for $l = 1, 2, 3$, where $A := -K_X$ is the generator of $\Cl (X) \cong \mathbb{Z}$ and $E$ is the exceptional divisor of $\varphi$.
By flitering off terms divisible by $x_j$ in the defining equation of $X$, we have
\[
(x_k^2 + \cdots) x_j = - x_k \bar{f} - \bar{g},
\]
where $\bar{f} = f (x_{i_1}, x_{i_2}, x_{i_3}, 0)$ and $\bar{g} = g (x_{i_1}, x_{i_2}, x_{i_3}, 0)$.
It follows that the section $x_j$ vanishes along $E$ to order at least $\deg f/a_k = (a_k+a_j)/a_k > a_j/a_k$.
This shows that $x_j$ lifts to a section of $-a_j K_Y$ and the claim is proved.

The variety $Z$ is explicitly described as
\[
Z = (y^2 + y f + x_j g = 0) \subset \overline{\mathbb{P}} := \mathbb{P} (a_{i_1},a_{i_2},a_{i_3},a_j,b),
\]
where $b := a_k + a_j$ and $x_{i_1},x_{i_2},x_{i_3},x_j,y$ are homogeneous coordinates of $\overline{\mathbb{P}}$ of weights $a_{i_1},a_{i_2},a_{i_3},a_j,b$, respectively.
The morphism $\psi \colon Y \to Z$ is induced by the map $\pi_X \colon X \dashrightarrow Z$ defined by
\[
(x_{i_1}\!:\!x_{i_2}\!:\!x_{i_3}\!:\!x_j\!:\!x_k) \mapsto
(x_{i_1}\!:\!x_{i_2}\!:\!x_{i_3}\!:\!x_j\!:\!x_k x_j).
\]
The biregular involution $\iota_Z$ of $Z$ which swaps 2 points in the $\pi$-fibers induces birational involutions $\iota_Y \colon Y \dashrightarrow Y$ and $\iota \colon X \dashrightarrow X$.
We obtain a diagram of the form \eqref{eq:QI}.

The inverse $\pi_X^{-1} \colon Z \dashrightarrow X$ is described as
\[
(x_{i_1}\!:\!x_{i_2}\!:\!x_{i_3}\!:\!x_j\!:\!y) \mapsto
(x_{i_1}\!:\!x_{i_2}\!:\!x_{i_3}\!:\!x_j\!:\!y/x_j) = (x_{i_1} x_j^{a_{i_1}/b}\!:\!x_{i_2} x_j^{a_{i_2}/b}\!:\!x_{i_3} x_j^{a_{i_3}/b}\!:\!x_j^{(a_j + b)/b}\!:\!y).
\]
Moreover, since $y/x_j = - g/(y+f)$ on $Z$, we have
\[
\begin{split}
(x_{i_1}\!:\!x_{i_2}\!:\!x_{i_3}\!:\!x_j\!:\!y/x_j) &=
(x_{i_1}\!:\!x_{i_2}\!:\!x_{i_3}\!:\!x_j\!:\!-g/(y+f)) \\
&= (x_{i_1} (y+f)^{a_{i_1}/b}\!:\!x_{i_2} (y+f)^{a_{i_2}/b}\!:\!x_{i_3} (y+f)^{a_{i_3}/b}\!:\!x_j (y+f)^{a_j/b}\!:\!-g).
\end{split}
\]
This shows that $\pi_X^{-1}$ is defined outside the set 
\[
(x_j = y = 0) \cap (y + f = g = 0) = (x_j = y = f = g = 0) =: \Xi_Z.
\]
Let $\Xi_Y$ be the strict transform of the set $(x_j = f = 0) \cap X = (x_j = f = g = 0) \subset X$ on $Y$.
Then $\psi (\Xi_Y) = \Xi_Z$ and $\dim \Xi_Z = \dim \Xi_Y - 1$.
It follows that the exceptional locus of $\psi$ is $\Xi_Y$.
By Proposition \ref{prop:QI}, the birational involution $\iota \colon X \dashrightarrow X$ is a Sarkisov self-link if and only if $\dim \Xi_Y = 1$.

If $x_j \mid f$, then $\dim \Xi_Y = 2$, that is, $\psi$ is divisorial, and hence $\varphi$ is not a maximal extraction by Proposition~\ref{prop:QI}.
This proves (1).

Suppose that $x_j \mid f$.
We show that $\dim \Xi_Y = 1$.
This is equivalent to showing that the subset 
\[
\overline{\Xi} := (\bar{f} = \bar{g} = 0) \subset \mathbb{A}^3 = \operatorname{Spec} (\mathbb{C} [x_{i_1},x_{i_2}, x_{i_3}])
\]
is of $1$-dimensional, where $\bar{f} = f (x_{i_1}, x_{i_2}, x_{i_3}, 0)$ and $\bar{g} = g (x_{i_1}, x_{i_2}, x_{i_3}, 0)$.
We have $\bar{f} \ne 0$ as a polynomial since $x_j \nmid f$.
If $\overline{\Xi}$ has an irreducible component which is of dimension $2$, then there is an irreducible component $f_1 = f_1 (x_{i_1}, x_{i_2}, x_{i_3})$ of $\bar{f}$ which divides $\bar{g}$.
In this case, we can write $f = f_1 f_2 + x_0 f_0$ and $g = f_1 g_1 + x_0 g_0$ for some homogeneous polynomials $f_1, f_2, g_1 \in \mathbb{C} [x_{i_1}, x_{i_2}, x_{i_3}]$ and $f_0, g_0 \in \mathbb{C} [x_{i_0}, x_{i_1}, x_{i_2}, x_j]$.
The polynomial defining $X$ can be written as
\[
x_k^2 x_j + x_k (f_1 f_2 + x_j f_0) + f_1 g_1 + x_j g_0 = x_j (x_k^2 + x_k f_0 + g_0) + f_1 (x_k f_2 + g_1).
\]
We see that $X$ is not quasismooth at any point of the nonempty set
\[
(x_j = x_k^2 + x_k f_0 + g_0 = f_1 = x_k f_2 + g_1 = 0) \subset X.
\]
This is absurd and thus $\dim \Xi_Y = 1$ and the assertion (2) follows from Proposition ~\ref{prop:QI}.
\end{proof}

By similar but more involved arguments, we can construct quadratic involutions for quasismooth Fano threefold WCIs of codimension $2$ and index $1$. 

\begin{Thm}[{\cite[Theorem~5.5]{okadaI}, \cite[Lemma~3.13]{KOWbirrigid}}]
Let $X = X_{d_1, d_2} \subset \mathbb{P} (a_{i_0}, a_{i_1}, a_{i_2}, a_{i_3}, a_j, a_k)$ be a quasismooth Fano threefold WCI of multidegree $(d_1, d_2)$ and of index $1$.
We assume that $x_k^2 x_{i_0}$ and $x_k x_{i_1}$ appear in the defining polynomials $F_1$ and $F_2$ of degree  $d_1$ and $d_2$ with nonzero coefficients, respectively.
Then, after replacing homogeneous coordinates, the defining equations can be written as
\[
\begin{split}
    F_1 &= x_k^2 x_{i_0} + x_k a + x_j^2 + b = 0, \\
    F_2 &= x_k x_{i_1} + x_j c + d,
\end{split}
\]
where $a, b, c, d \in \mathbb{C} [x_{i_0}, x_{i_1}, x_{i_2}, x_{i_3}]$ are homogeneous polynomials.
Let $\mathbf{p} \in X$ be the point at which only the homogeneous coordinate $x_k$ does not vanish and let $\varphi \colon Y \to X$ be the Kawamata blowup at $\mathbf{p}$.
\begin{enumerate}
\item If $c = 0$ as a polynomial, then $\varphi$ contracts a divisor and $\mathbf{p}$ is not a maximal extraction.
\item If $c \ne 0$, then there exists a birational involution $\varphi \colon X \dashrightarrow X$ which is a Sarkisov self-link.
\end{enumerate}
\end{Thm}

We explain the construction of elliptic involutions.
Let $U$ be a normal projective $\mathbb{Q}$-factorial variety with only terminal singularities.
Assume that $-K_U$ is semiample and the anticanonical morphism $\tau \colon U \to S$ is of relative dimension $1$.
We assume in addition that there is a prime divisor $F$ on $U$ which is a rational section of $\tau$.
Then the generic fibre $U_{\eta}$ of $\tau$ is an elliptic curve with a $K (S)$-rational point $\mathbf{p}_F$ corresponding to the rational section $F$.
The reflection of $U_{\eta}$ with respect to the point $\mathbf{p}_{F}$ gives rise to a birational involution $\iota_U \colon U \dashrightarrow U$ over $S$.
The map $\iota_U$ is an isomorphism in codimension $1$ since $U$ has only terminal singularities and $K_U$ is nef over $S$.
We call $\iota_U$ the elliptic involution with respect to $F$.

\begin{Prop}[Elliptic involution]
\label{prop:EI}
Let $X$ be a Fano threefold with $\rho (X) = 1$ and let $\varphi \colon Y \to X$ be a divisorial extraction.
Suppose that there exist a birational morphism $\psi \colon W \to Y$ from a normal projective variety $W$, a birational map $\chi \colon W \dashrightarrow U$ to a normal projective $\mathbb{Q}$-factorial variety $U$ with only terminal singularities which is an isomorphism in codimension $1$ and a prime divisor $F$ on $W$ which is contracted by $\psi$ satisfying the following properties.
\begin{enumerate}
    \item The anticanonical divisor $-K_U$ is semiample and the anticanonical morphism $\tau \colon U \to S$ is of relative dimension $1$.
    Moreover the divisor $\chi_*F$ is a rational section of $\tau$.
    \item Any prime divisor on $W$ other than $F$ which is contracted by $\psi$ is contracted by $\tau \circ \chi$.
\end{enumerate}
\begin{equation} \label{eq:EI}
\vcenter{
\xymatrix{
& \ar[ld]_{\psi} W \ar@{-->}[r]^{\chi} & U \ar[rdd]_{\tau} \ar@{-->}[rr]^{\iota_U} & & \ar[ldd]^{\tau} U & \ar@{-->}[l]^{\chi} W \ar[rd]^{\psi} & \\
Y \ar[d]_{\varphi} & & & & & & Y \ar[d]^{\varphi} \\
X & & & S & & & X
}
}
\end{equation}
Then the elliptic involution $\iota_U \colon U \dashrightarrow U$ with respect to $\chi_*F$ induces a birational map $\iota_Y \colon Y \dashrightarrow Y$.
If $\varphi$ is a maximal extraction and $\iota_Y$ is not an isomorphism, then the birational map $\iota := \varphi \circ \iota_Y \circ \varphi^{-1} \colon X \dashrightarrow X$ is a Sarkisov self-link.
\end{Prop}

\begin{proof}
By the property (1) and the preceding arguments, there is a birational involution $\iota_U \colon U \dashrightarrow U$ which is an isomorphism in codimension $1$.
By the assumption that $\chi$ is an isomorphism in codimension $1$ and by the property (2), the birational involution $\iota_Y$ of $Y$ induced by $\iota_U$ is also an isomorphism in codimension $1$.
The rest follows from \cite[Lemma~2.26]{okadaII}.
\end{proof}

We give an example of elliptic involutions for Fano threefold weighted hypersurfaces.
Elliptic involutions can also be seen as quadratic involutions: we refer readers to \cite[\S 4.2]{CPR} and \cite[\S~5.2]{okadaI}, where elliptic involutions for Fano threefold weighted hypersurfaces and WCIs of codimension 2, respectively, are constructed as quadratic involutions.
There are more involved elliptic involutions which are called invisible involutions (see \cite[\S~4.3]{chel} and \cite[\S~7]{okadaII}).

\begin{Ex}
Let $X = X_{18} \subset \mathbb{P} (1, 1, 4, 6, 7)$ be a quasismooth weighted hypersurface of degree $18$ and let $\mathbf{p} = (0\!:\!0\!:\!1\!:\!0\!:\!0) \in X$ be the singular point of type $\frac{1}{4} (1, 1, 3)$.
By a suitable choice of homogeneous coordinates $x_0, x_1, x_2, x_3, x_4$ of weights $1, 1, 4, 6, 7$, respectively, we can write the equation of $X$ as
\[
x_2^3 x_3 + x_2 (x_4^2 + f_{14} (x_0,x_1,x_3)) + x_4 g_{11} (x_0, x_1, x_3) + h_{18} (x_0, x_1, x_3) = 0, 
\]
where $f_{14}, g_{11}, h_{18} \in \mathbb{C} [x_0, x_1, x_3]$ are homogeneous polynomials of degree $14, 11, 18$, respectively.
Let $X \dashrightarrow \mathbb{P} (1, 1, 6) =: S$ be the projection to the coordinates $x_0, x_1, x_3$.
Then the the point $\mathbf{q}$ is the unique indeterminacy locus of the induced map $Y \dashrightarrow S$, where $\mathbf{q} \in Y$ is the point which is mapped isomorphically to the point $(0\!:\!0\!:\!0\!:\!0\!:\!1) \in X$ of type $\frac{1}{7} (1, 1, 6)$.
Let $\psi \colon W \to Y$ be the Kawamata blowup of $Y$ at $\mathbf{q}$.
Then $\psi$ resolves the indeterminacy of the map $Y \dashrightarrow S$.
It can be seen that the induced morphism $\eta \colon W \to S$ is the anticanonical morphism (by observing that the sections $x_0, x_1, x_3$ lift to (pluri-)anticanonical sections on $W$) and that the $\psi$-exceptional divisor $F \cong \mathbb{P} (1, 1, 6)$ is the section of $\eta \colon W \to S$.
Thus, we obtain the diagram of the form \eqref{eq:EI} with $U = W$ and $\chi$ being the identity.
By Proposition \ref{prop:EI}, there are associated birational involutions $\iota_Y \colon Y \dashrightarrow Y$ and $\iota \colon X \dashrightarrow X$. 
Although it is not easy to see this, we can show that the birational involution $\iota_Y$ is not biregular if $\mathbf{p}$ is a maximal centre (see \cite[Proof of Theorem 4.2.6]{chel}).  
\end{Ex}

As far as the authors know, every Sarkisov self-link that appears in the study of birational rigidity/solidity of a Fano variety of Picard rank $1$ is a birational involution and we ask the following.

\begin{Quest}
Is there a Sarkisov self-link $\iota \colon X \dashrightarrow X$ of a Fano variety $X$ with $\rho (X) = 1$ which is not a birational involution?
\end{Quest}

\subsection{Quartic threefolds}

As a demonstration of the method of maximal singularities, we give a complete proof of the following foundational result of birational (super-)rigidity. 

\begin{Thm}[\cite{iskmanin}]
\label{thm:IskManin}
Any smooth quartic threefold is birationally superrigid.
\end{Thm}

\begin{proof}
Let $X$ be a smooth quartic threefold.
We set $A = -K_X$ which is the class of a hyperplane section on $X$.
We assume that $X$ is not birationally superrigid.
By Theorem~\ref{thm:NFI}, there exists a maximal centre $\Gamma \subset X$ with respect to some mobile linear system $\mathcal{M} \sim n A$.
The center $\Gamma$ is either a point or a curve.

Suppose that $\Gamma$ is a point $\mathbf{p} \in X$.
Let $D_1, D_2 \in \mathcal{M}$ be general members.
Then, for the effective $1$-cycle $D_1 \cdot D_2$, we have $\mult_{\mathbf{p}} (D_1 \cdot D_2) > 4 n^2$ by Theorem \ref{thm:4nineq}.
We can take a hyperplane $H \in |A|$ such that $H \cap D_1 \cap D_2$ is a finite set of points.
Then we have
\[
4 n^2 = (H \cdot D_1 \cdot D_2) \ge \mult_{\mathbf{p}} (H) \mult_{\mathbf{p}} (D_1 \cdot D_2) > 4 n^2,
\]
where the first inequality follows from \cite[Corollary 12.4]{FultonIT}.
This is impossible.

It follows that $\Gamma$ is an irreducible and reduced curve $C$.
By Lemma \ref{lem:exclCurvebydeg}, we have $\deg C \leq 3$ and $C$ is one of the following (cf.\ \cite[IV, Exercises~3.4]{Hart}).
\begin{itemize}
\item $C$ is a line.
\item $C$ is a conic.
\item $C$ is a plane cubic curve.
\item $C$ is a twisted cubic curve.
\end{itemize}

We consider the case where $C$ is not a plane cubic curve.
Then $C$ is a smooth rational curve and the blowup $\varphi \colon Y \to X$ along $C$ is the unique divisorial extraction centered at $C$ (cf.\ \cite[Proposition~4.1.3]{kawakitabook}) and the pair $(X, \frac{1}{n} \mathcal{M})$ is not canonical at its exceptional divisor $E$.
We have
\[
\begin{split}
(\varphi^*A)^3 &= A^3 = 4, \\
(\varphi^*A)^2 \cdot E &= 0, \\ 
(\varphi^*A) \cdot E^2 &= - \deg C, \\ 
E^3 &= - \deg \mathcal{N}_{C/X} = - \deg C + 2.
\end{split}
\]
We set $B := -K_Y = \varphi^*A - E$.
If $C$ is a line, then $|\mathcal{I}_C (A)|$ isolates $C$, and hence $M := \varphi^*A - E$ is nef.
If $C$ is either a conic or a twisted cubic curve, then $|\mathcal{I}_C (2 A)|$ isolates $C$ since $C$ is cut by quadrics, and hence $M := 2 \varphi^*A - E$ is nef.
We compute
\[
M \cdot B^2 = 
\begin{cases}
0, & \text{if $C$ is either a line or a conic}, \\
-1, & \text{if $C$ is a twisted cubic}.
\end{cases}
\]
By Lemma \ref{lem:exclNE}, $C$ is not a maximal center.

It remains to consider the case where $C$ is a plane cubic curve.
We set $\mathcal{H} := |\mathcal{I}_C (A)|$.
We see that $\mathrm{Bs} (\mathcal{H}) = C$.
Let $S, T \in \mathcal{H}$ be general members.
Then $S$ and $T$ are smooth surfaces and we have $S|_T = C + \ell$, where $\ell$ is a line.
Clearly we have $C \cdot \ell \ge A \cdot \ell = 1$.
By Lemma \ref{lem:exclCurveint}, $C$ is not a maximal center.
We have derived a contradiction.
Therefore $X$ is birationally superrigid.
\end{proof}

\section{Constructing Sarkisov links}
\label{sec:SarkisovLinks}

In this section we discuss an approach to construct Sarkisov links from a Fano threefold, necessarily of Picard rank $1$. For simplicity we focus on the case where the centre of the divisorial extraction is a point.

We start by briefly recalling several cones in birational geometry. Let $Y$ be a projective variety. We denote by $\mathrm{N}^1(Y)$ the (finitely generated) abelian group of Cartier divisors modulo numerical equivalence and we call it the Néron-Severy space of $Y$.  Then $\mathrm{N}^1(Y)_{\mathbb R} :=\mathrm{N}^1(Y)\otimes_{\mathbb Z} \mathbb R$ is an $\mathbb R$-vector space. The cone generated by nef divisors in $\mathrm{N}^1(Y)_{\mathbb R}$ is denoted $\mathrm{Nef}(Y)$. Similarly, the closure of the cone of effective divisors (resp., movable divisors) is denoted $\overline{\mathrm{Eff}}(Y)$
(resp., $\overline{\mathrm{Mov}}(Y)$). By a small $\mathbb Q$-factorial modification, or SQM,  of a normal projective variety $Y$ we mean a small (i.e., isomorphic in codimension one) birational
map $Y \dashrightarrow Y^+$ to another normal, $\mathbb Q$-factorial projective variety $Y^+$. Notice that an SQM can be an isomorphism.

We recal the definition of a Mori dream space. See \cite{CastravetMDS} for a survey around this topic.

\begin{Def}[{\cite[Definition~1.10]{mdsGIT}}]
A normal projective variety $Y$ is a \textbf{Mori dream space} (MDS for short) if the following hold
\begin{itemize}
	\item $Y$ is $\mathbb{Q}$-factorial and $\mathrm{Pic}(Y)_{\mathbb{Q}} = N^1(Y)$;
	\item $\mathrm{Nef}(Y)$ is the affine hull of finitely many semi-ample line bundles;
	\item There exists a finite collection of small $\mathbb{Q}$-factorial modifications $\tau_i \colon Y \dashrightarrow Y_i$ such that each $Y_i$ satisfies the previous points and $\overline{\mathrm{Mov}}(Y)= \bigcup \tau_i^*(\mathrm{Nef}(Y_i))$.
\end{itemize}
\end{Def}

\begin{Rem}
    $\mathbb Q$-factorial varieties of Fano type are Mori dream spaces by \cite[Cor~1.3.1]{BCHM}. Examples include toric varieties. 
\end{Rem}

First observe that if $\varphi \colon Y \rightarrow X$ initiates a Sarkisov link then $Y$ is a Mori dream space, see for instance, \cite[Lemma~2.9]{hamidquartic}. On the other hand, by \cite[Proposition~2.11]{mdsGIT}, the birational contractions of a Mori dream space $Y$ are induced from toric geometry. In particular, there is an embedding $Y \subset T$ into a quasismooth projective toric variety for which every Mori chamber of $Y$ is a finite union of Mori chambers of $T$. 

We briefly explain the idea of the embedding $Y \subset T$. Let $X \subset \mathbb{P}$ be a quasismooth Fano 3-fold in a weighted projective space that we denote by $\mathbb{P}$ for simplicity. If $\mathbf{p} \in X \subset \mathbb{P}$ is a germ of a point, we perform a toric blowup  $\Phi \colon T \rightarrow \mathbb{P}$ centred at $\mathbf{p}$ with the constraint that it restricts locally around $\mathbf{p}$ to a terminal divisorial extraction from $\mathbf{p} \in X$ as explained in Section \ref{sec:termextr},
\[
\varphi \colon E \subset Y \rightarrow \mathbf{p} \in X, 
\]
where $Y:= \Phi_*^{-1}X$ and $\varphi$ is the restriction $\Phi|_{\Phi_*^{-1}X}$. This is summarised in the following diagram
\[
\begin{tikzcd}[column sep = 4em]
E \arrow[swap]{d}{\varphi|_E} \arrow[swap,hookrightarrow]{r}{} & Y \arrow[swap]{d}{\varphi} \arrow[swap,hookrightarrow]{r}{} \arrow[dashed]{r}{} & T \arrow[swap]{d}{\Phi}    \\
 \mathbf{p} \arrow[swap,hookrightarrow]{r}{} & X \arrow[hookrightarrow]{r}{} & \mathbb{P}   
\end{tikzcd}
\]

Once we have a candidate for $T$, we describe its birational geometry. That is, we look at the closure of the movable and effective cones of $T$, denoted $\overline{\Mov}(T)$ and $\overline{\Eff}(T)$, respectively and how these are related. Notice that since $\Phi$ is a divisorial extraction from $\mathbb{P}$, the movable cone of $T$ is strictly contained in the effective cone of $T$. The chambers $f_i^*(\Nef(T_i))$ and their faces form a fan supported at $\overline{\Mov}(T)$. The number of cones in this fan is finite by \cite[Proposition~1.11~(2)]{mdsGIT} and are in one-to-one correspondence with contracting rational maps $g \colon T \rat Q$, with $Q$ normal and projective via
\[
g \colon T \rat Q  \quad \leftrightsquigarrow \quad g^*(\Nef(Q)) \subset \overline{\Mov}(T)
\]
by \cite[Proposition~1.11~(3)]{mdsGIT}. In particular, if $T_i$ and $T_{i+1}$ are the ample models (see \cite[Definition~3.6.5]{BCHM}) of adjacent chambers, then they are related by a small $\mathbb{Q}$-factorial modification that we denote by $\tau_i$. Since in our case $\Pic(T)=2$, we can think of these cones in $\mathbb{R}^2$ and hence there is a completely determined sequence of birational transformations 
\[
\begin{tikzcd} 
& T \arrow[swap]{dl}{\Phi} \arrow[swap]{dr}{f} \arrow[dashed]{rr}{\tau} & &  T_1\arrow[swap]{dr}{f_1} \arrow{dl}{g} \arrow[dashed]{r}{\tau_1} & \cdots \arrow[dashed]{r}{\tau_n}   & T'\arrow{dr}{\Phi'}\arrow{dl}{g_{n-1}}  & \\
 \mathbb{P} &  & \mathcal{F} &  & \cdots  &  & \mathcal{F}' 
\end{tikzcd}
\]
starting from $\Phi \colon T \rightarrow \mathbb{P}$ called the \textbf{2-ray game}. Notice that we can always carry a 2-ray game in a $\mathbb Q$-factorial toric variety of Picard rank 2. See also \cite[Definition~4.2.5]{kawakitabook}. As explained, these maps are induced from a chamber decomposition of $\overline{\Eff}(T)$
There are two possible outcomes of running the 2-ray game on $T$. 
\begin{enumerate}
	\item The closure of the movable cone of $T$ is stricly contained in the interior of the cone of pseudo-effective divisors of $T$. In this case, the map $\Phi'$ is birational and is, in fact, a divisorial contraction to a closed subvariety. 
	\item The closure of the movable cone of $T$ is \emph{not} stricly contained in the interior of the cone of pseudo-effective divisors of $T$. In this case, there is a divisor class $D$ in the boundary of $\overline{\Mov}(T)$ which is not big. Hence the map $\Phi'$ given by the positive multiples $|mD|$ is not birational and $\Phi'$ is a fibration.
\end{enumerate}

It is then natural to relate the effective cone decomposition of $T$ with that of $Y$. It turns out, that in general, these do not coincide. For example the $\mathrm{Nef}(T)$ can be strictly smaller than $\mathrm{Nef}(Y)$. See for instance \cite{hassett2002weak} and \cite[Lemma~3.6]{tiagoliviaerik}.

\begin{Def}[{\cite[Definition~3.5]{hamidmds}}] \label{def:follow}
Let $Y \subset T$ be $\mathbb{Q}$-factorial projective varieties and suppose $\Pic(T)=2$. We say that \textbf{$Y$ 2-ray follows $T$} or that \textbf{$Y$ follows $T$} if the 2-ray game on $T$ restricts to a 2-ray game on $Y$ where $Y_i \subset T_i$ is given by the ideal $I_Y$ and the small $\mathbb{Q}$-factorial modification $T_i \rat T_{i+1}$ over $\mathcal{F}_i$ restricts to a small $\mathbb{Q}$-factorial modification $Y_i \rat Y_{i+1}$ over $\mathcal{F}_i|_{Y_i}$. 
\end{Def}

\begin{Ex} \label{ex:quarticsarkisovlink}
    We consider terminal $\mathbb Q$-factorial quartic threefolds of the form 
    \begin{equation}
\label{eq:quartic}    X\colon (x_0^2x_1x_2+x_0f_3+f_4=0) \subset \mathbb P^4
    \end{equation}
    where $f_3$ and $f_4$ are homogeneous polynomials in $\mathbb C[x_1,x_2,x_3,x_4]$ of degrees $3$ and $4$, respectively. Either $(f_3,f_4)$ is a regular sequence in $\mathbb C[x_1,x_2,x_3,x_4]$ or there exists a non-constant homogeneous polynomial  $h_i \in \mathbb C[x_1,x_2,x_3,x_4]$ of degree $1 \leq i\leq 4$ such that $h_i$ is a common factor of both $f_3$ and $f_4$ and $\mathrm{gcd}\big(\frac{f_3}{h_i},\frac{f_4}{h_i}\big)=1$. 
    
    We claim that if $(f_3,f_4)$ does not form a regular sequence, then $f_3$ is identically zero. Indeed, if $(f_3,f_4)$ is not a regular sequence, then we can write $X$ as
         $$
    X\colon (x_0^2x_1x_2+h_i(x_0f_{3-i}+f_{4-i})=0) \subset \mathbb P^4.
    $$
    Let $g:=x_0f_{3-i}+f_{4-i}$. The Jacobian of $X$ is 
    $$
    \mathrm{J}(X)=\begin{pmatrix}
 2 x_0x_1x_2+h_i\cdot f_{3-i}\\ x_0^2x_2 +g\cdot \partial_{x_1} h_i + h_i\cdot \partial_{x_1}g \\
 x_0^2x_1 +g\cdot \partial_{x_2} h_i + h_i\cdot \partial_{x_2}g \\   g\cdot \partial_{x_3} h_i + h_i\cdot \partial_{x_3}g \\ g\cdot \partial_{x_4} h_i + h_i\cdot \partial_{x_4}g
\end{pmatrix}
    $$
    where the $i^{\text{th}}$ row corresponds to the partial derivative of the defining equation of $X$ with respect to the variable $x_{i-1}$.  Notice that $\partial_{x_i} g = x_0\cdot \partial_{x_i} f_{3-i}+\partial_{x_i}f_{4-i}$ for $i\geq 1$ and thus, if $f_{4-i}$ is non-constant, we can readily see that $X$ is singular along the locus $(x_0=h_i=g=0)$ which is at least $1$-dimensional, contradicting the fact that $X$ is terminal. Hence $f_{4-i}\in \mathbb C^*$ which implies that $\mathrm{deg} (h_i)=4$ and so $f_3\equiv 0$, as we claimed.

    Hence, terminal $\mathbb Q$-factorial quartic threefolds of the form \ref{eq:quartic} are the fibres of the family flat $\pi \colon \mathcal X \rightarrow \mathbb C$
$$
\mathcal{X} \colon (x_0^2x_1x_2+t\cdot x_0f_3+f_4=0) \subset \mathbb P^4\times \mathbb C_t
$$
where $(f_3,f_4)$ form a regular sequence.
We let $X_t:=\pi^{-1}(t)$ be the fibre over the point $t \in \mathbb C$. If $t$ is the general point of $\mathbb C$ we denote $\pi^{-1}(t)$ by $X$ and if $t=0$ we denote the central fibre $\pi^{-1}(0)$ by $X_0$. Notice that when $f_3$ and $f_4$ are general then $X$ is the quartic threefold from \cite{cortimella}.

We construct Sarkisov links centred at the point $\mathbf{p_{x_0}}=[1:0:0:0:0]\in X_t$. Suppose $X$ has a $cA_2$ singularity at $\mathbf{p_{x_0}}$ which is locally analytically equivalent to 
$$
(0 \in xy+z^3+t^3=0) \subset \mathbb C^4.
$$

We choose (global) blowup weights which restrict to a divisorial contraction of the germ $(0 \in xy+z^3+t^3=0) \subset \mathbb C^4$. By Kawakita, Theorem \ref{thm:divcontcA}, we have $\mathrm{wt} (x,y,z,t) = (2,1,1,1)$ or $(1,2,1,1)$. So the choice $\mathrm{wt}(x_1,x_2,x_3,x_4)=(2,1,1,1)$ is natural. That is, we consider the weighted blowup $\sigma \colon T \rightarrow \mathbb P^4$ centred at $\mathbf{\mathbf{p_{x_0}}}$ with weights $\mathrm{wt}(x_1,x_2,x_3,x_4)=(2,1,1,1)$ and exceptional divisor $E\simeq \mathbb P(2,1,1,1)$. Let $H:=\sigma^*\mathcal{O}_{\mathbb P^4}(1)$ denote the pullback of a hyperplane section from $\mathbb P^4$. Moreover, let $D_{x_1}$ be the strict transform of $(x_1=0)$ in $T$. Then, the morphism $\sigma'$ induced by the linear system $|H-E|$ is a divisorial contraction $\sigma'\colon T\rightarrow \mathbb P(2,1,1,1)$ that contracts the divisor $D_{x_1}$ to a projective plane.     Let $Y$ be the strict transform of $X\subset \mathbb P^4$ under $\sigma$. In other words, the threefold $Y$ is a weighted blowup of $X$ at the point $\mathbf{p_{x_0}}$ with weights $\mathrm{wt}(x_1,x_2,x_3,x_4)=(2,1,1,1)$. Then $Y \in |4H-3E|$ is given by
    $$
    Y \colon (x_0^2x_1x_2+t\cdot x_0f_3(x_1u,x_2,x_3,x_4)+uf_4(x_1u,x_2,x_3,x_4)=0) \subset T,
    $$
    where $u$ is the section correcponding to the exceptional divisor $E$, and $\sigma'{\vert_Y}$ is a contraction to
$$
Z \colon (x_0^2x_2+t\cdot x_0f_3(u,x_2,x_3,x_4)+uf_4(u,x_2,x_3,x_4)=0) \subset \mathbb P(2,1,1,1,1)
$$
which is not $\mathbb Q$-factorial for any $t \in \mathbb C$. Hence, the threefold $Y$ does not follow $T$ in the sense of Definition \ref{def:follow}. One way to overcome this problem is to enlarge the toric variety $T$ by adding effective divisors corresponding to sections of $Y$. We define 
$$
\eta :=\frac{x_0x_1x_2+t\cdot f_3(x_1u,x_2,x_3,x_4)}{u}=-\frac{f_4(x_1u,x_2,x_3,x_4)}{x_0} \in H^0(Y,\mathcal{O}_{Y}(3H-4E)).
$$
Notice how $\eta$ is indeed a regular section of $Y$. Define the new toric variety $T^u$ as the quotient (see \cite[Theorem~5.1.11]{Coxbook}) 
$$
\frac{\mathrm{Spec} \big(\mathbb C[u,x_0,x_1,x_2,x_3,x_4,\eta] \big) \setminus \mathcal{Z}(I_{\text{irr}})}{\mathbb C^*\times \mathbb C^*} 
$$
where the irrelevant ideal is $I_{\text{irr}}=(u,x_0)\cap(x_1,x_2,x_3,x_4,\eta)$ and the $(\mathbb C^*)^2$ action is
$$
(\lambda,\mu)\cdot (u,x_0,x_1,x_2,x_3,x_4,\eta) = (\mu^{-1}u,\lambda^{}x_0,\lambda^{}\mu^{2}x_1,\lambda^{}\mu^{}x_2,\lambda^{}\mu^{}x_3,\lambda^{}\mu^{}x_4,\lambda^{3}\mu^{4}\eta). 
$$
In fact, $T^u$ is the weighted blowup of $\mathbb P(1^5,3)$ at the point $[1:0:\cdots : 0]$ with coordinate weights $(2,1,1,1,4)$. The decomposition of the effective cone of $T^u$ into Mori chambers together with the corresponding birational transformations is given in Figure \ref{fig:test1} and Figure \ref{fig:test2}, respectively, see \cite[Chapter~15]{Coxbook} for an extensive treatment. The map
$$f \colon T^u \rightarrow \mathcal F=\mathrm{Proj}\big(\bigoplus_{m\geq 0}H^0(T^u,m(H-E))\big)
$$
is given in coordinates by $(x_2:x_3:x_4:ux_1:x_0x_1:u\eta:x_0\eta)$. Therefore the image of $f$ is 
$$
\mathcal F \simeq (v_4v_7-v_5v_6=0) \subset \mathrm{Proj} \big(\mathbb C[v_1,v_2,v_3,v_4,v_5,v_6,v_7]\big) \simeq  \mathbb P(1,1,1,1,2,3,4).
$$

It follows that $f$ contracts the locus $L\colon (\eta=x_1=0)\subset T^u$, which is isomorphic to the blowup of $\mathbb P^3$ at a point, and is an isomorphism otherwise. Notice that $L$ admits a fibration to $\mathbb P^2$ and so $f$ is a \emph{small} contraction  which contracts $L$ to the base of its fibration. Similarly we have,
$$
T_1^u:=\frac{\mathrm{Spec} \big(\mathbb C[u,x_0,x_1,x_2,x_3,x_4,\eta] \big) \setminus \mathcal{Z}(I_{\text{irr}})}{\mathbb C^*\times \mathbb C^*}
$$
where the irrelevant ideal is $I_{\text{irr}}=(u,x_0,x_2,x_3,x_4)\cap(\eta,x_1)$ and the $(\mathbb C^*)^2$ is the same as above. In the same way, $T_1^u$ admits a small contraction $g \colon T_1^u \rightarrow \mathcal F=\mathrm{Proj}\big(\bigoplus_{m\geq 0}H^0(T_1^u,m(H-E))\big) 
$
which contracts $L_1\colon (u=x_0=0)$ to $\mathbb P^2$ and is an isomorphism otherwise. The map $\tau:=g\circ f^{-1}$ is a small $\mathbb Q$-factorial modification of $T^u$, in fact, the only one besides the identity, see \cite{mdsGIT}, and it swaps the fibres over $\mathbb P^2$ of $L$ to those of $L_1$. Finally $\sigma'\colon T_1^u\rightarrow \mathbb P(1^4,2^2)$ is a divisorial contraction of $(x_1=0)$ to a point of type $\frac{1}{2}(1,1,1,1,0)$.

    \begin{figure}[h!]
\centering
\begin{minipage}{.45\textwidth}
  \centering
  \begin{tikzpicture}[scale=2]
  \coordinate (A) at (0, 0);
  \coordinate [label={left:$E$}] (E) at (0, -0.7);
  \coordinate [label={right:$H$}] (K) at (1, 0);
\coordinate [label={right:$H-E$}] (y) at (1,1);
\coordinate [label={right:$3H-4E$}] (M2) at (0.5,1.2);
\coordinate [label={above:$H-2E$}] (E') at (0.15,1.2);
\draw[fill=gray!30,draw=none]    (0,0) -- ++(1,0) -- ++(0,1) ;
    \draw [very thick] (A) -- (E);
    \draw [very thick] (A) -- (y);
    \draw [very thick] (A) -- (E');
    \draw [very thick] (A) -- (K);
    \draw (0.7,0.35) node{ \tiny \tiny$\mathrm{Nef}(T^u)$} ;
   \draw [very thick] (A) -- (M2);
\end{tikzpicture}
  \captionof{figure}{}
  \label{fig:test1}
\end{minipage}%
\hfill
\begin{minipage}{.5\textwidth}
  \centering
  \[
        \begin{tikzcd}[ampersand replacement=\&,column sep = 1.5em]
             \&   T^u \ar[dl,swap, "\sigma"] \ar[rr,dashed,"\tau"] \ar[rd,swap,"f"]  \& \& T^u_1 \ar[rd,"\sigma'"] \ar[ld,"g"]  \&  \\
             \mathbb{P}(1^4,3)  \& \& \mathcal{F} \&  \& \mathbb{P}(1^4,2^2)
        \end{tikzcd}
    \]
  \captionof{figure}{}
  \label{fig:test2}
\end{minipage}
\end{figure}

The projection $\pi_T \colon T^u \dashrightarrow T$ given by 
$$
(u,x_0,x_1,x_2,x_3,x_4,\eta)\mapsto (u,x_0,x_1,x_2,x_3,x_4)
$$
has base locus $(x_1=x_2=x_3=x_4=0)$ and its image is $T$. The inverse of $\pi$ is called an \emph{unprojection} and the term was coined by Miles Reid. We now consider the complete intersection
$$
  Y^u\colon 
  \begin{cases}
    \eta u -x_0x_1x_2-t\cdot f_3(x_1u,x_2,x_3,x_4=0 \\
    \eta x_0+f_4(x_1u,x_2,x_3,x_4)=0
  \end{cases}
   \subset T^u
$$
The point is that the map $\pi_T$ restricts to a morphism on $Y^u$ which is actually an isomorphism onto its image, $Y$. We claim that if we restrict the maps from $T^u$ to $Y^u$ we get a Sarkisov link. Recall that $f$ contracts the locus $L\subset T^u$ and is an isomorphism otherwise. Moreover $Y^u\cap L$ restricts to a $\mathbb P^1$-bundle over the complete intersection of a cubic and a quartic curves in $\mathbb P^2$ given by $\{f_3(0,x_2,x_3,x_4)=f_4(0,x_2,x_3,x_4)=0$. Hence, if $(x_1,f_3,f_4)$ form a regular sequence these are $12$ coplanar points. Otherwise these two curves share a component which is a curve and $L$ is a non Cartier divisor in $Y^u$, which contradicts $\mathbb Q$-factoriality of $Y$. We therefore assume from now on that $(x_1,f_3,f_4)$ forms a regular sequence. Let $Y_1^u \subset T^u_1$ be the threefold given by the same equations of $Y^u$. In the same way, $Y_1^u\cap L_1$ is a $\mathbb P^1$-bundle over the same $12$ points (more precisely a flopping contraction). In other words $\tau\vert_{Y^u}$ swaps $12$ rational curves $C_i$ with other $12$, $C_i'$, and $K_{Y_1}\cdot C_i = K_{Y_1^u}\cdot C_i'=0$, so $\tau$ restricts to $12$ simultaneous flops.

Finally $Y_1^u\cap(x_1=0)\simeq \mathbb P^2$ and $\sigma'\vert_{Y_1^u}$ is the Kawamata blowup of 
$$
  Z_{3,4}\colon 
  \begin{cases}
    \eta u -x_0x_2-t\cdot f_3(u,x_2,x_3,x_4)=0 \\
    \eta x_0+f_4(u,x_2,x_3,x_4)=0
  \end{cases}
   \subset \mathbb P(1^4,2^2)
$$
along the point $\mathbf{p_{\eta}}\sim\frac{1}{2}(1,1,1)$. This is the unique terminal extraction from $\mathbf{p_{\eta}}$ by Theorem \ref{thm:kwbl}. The Sarkisov link is then

  \[
        \begin{tikzcd}[ampersand replacement=\&,column sep = 2.5em]
             \&   Y^u \ar[dl,swap, "{(2,1,1,1)}"] \ar[rr,dashed,"\text{12 flops}"]   \& \& Y^u_1 \ar[rd,"{\frac{1}{2}(1,1,1)}"]  \&  \\
            \mathbf{p_{x_0}} \in X  \& \& \&  \& Z_{3,4} \ni \mathbf{p_{\eta}}
        \end{tikzcd}
    \]

    We finish this example by giving a Sarkisov link from the central fibre of $\pi$. Recall that, in this case, $X_0$ is a quartic threefold with a $cA_3$ singularity provided that $f_4$ is general. We do a weighted blowup of $\mathbb P^4$ at the coordinate point $\mathbf{p_{x_0}}$ with weights $\mathrm{wt}(x_1,x_2,x_3,x_4)=(2,2,1,1)$. Let $\sigma \colon T\rightarrow \mathbb P^4$ be such a blowup. Then $\sigma$ restricts to a $(2,2,1,1)$-weighted blowup of $X$ along $\mathbf{p_{x_0}}$, which is a terminal extraction with discrepancy $1$ by Theorem \ref{thm:divcontcA} and it initiates a Sarkisov link where $4$ rational curves on $Y$ are flopped to a del Pezzo fibration of degree $2$, as depicted in the following diagram.
      \[
        \begin{tikzcd}[ampersand replacement=\&,column sep = 2.5em]
             \&   Y_0 \ar[dl,swap, "{(2,2,1,1)}"] \ar[rr,dashed,"\text{4 flops}"]   \& \& {Y_0}_1 \ar[d,"\sigma'"]  \&  \\
            \mathbf{p_{x_0}} \in X_0  \& \& \& \mathbb P^1 \& 
        \end{tikzcd}
    \]
\end{Ex}

\section{The big picture on Fano threefold WCI}
\label{sec:bigpicture}

\subsection{Fano index equal to 1}

After Reid's discovery of the 95 families of $K3$ surfaces with canonical singularities in $\mathbb P(a_1,a_2,a_3,a_4)$ \cite{reid1980canonical}, Fletcher (a PhD student of Reid at the time) found that there is a bijection between this set and the set of families of quasismooth Fano threefold hypersurfaces in $\mathbb P(1,a_1,a_2,a_3,a_4)$ of Fano index $1$ and with terminal singularities \cite[Lemma~16.4]{CPR}, where the latter set exhausts all the quasismooth Fano threefold weighted hypersurfaces of Fano index $1$ (see \cite{JKWeightedHyp} and \cite{CCC}). We refer to these, and other subsequent lists, as Reid-Fletcher lists. Among the 95, there are 2 families of smooth weighted hypersurfaces of index $1$: the family of smooth quartic threefolds in $\mathbb{P}^4$ and the family of smooth sextic hypersurfaces in $\mathbb{P} (1,1,1,1,3)$ which are also known as sextic double solids. Birational superrigidity of smooth Fano threefolds in these 2 families is proved by Iskovskikh-Manin \cite{iskmanin} (see also Theorem~\ref{thm:IskManin}) for quartics and by Iskovskikh \cite{iskovskikh1979birational} for sextic double solids.

Corti, Pukhlikov and Reid \cite{CPR} studied birational geometry of the remaining 93 families systematically and obtained the following

\begin{theorem}[{{\cite[Theorem~1.3, Theorem~3.2]{CPR}}}]
    Let $X$ be a quasismooth Fano threefold weighted hypersurface with $\iota_X=1$. Suppose that $X$ is general in its family. Then $X$ is birationally rigid.
\end{theorem}

However complete the result was, the authors still acknowledged that the task was not over and conjectured that the same should hold for \emph{any} quasismooth member of the deformation family of $X$. The conjecture was proven only recently in \cite{chel}:

\begin{theorem}[{{\cite[Main Theorem]{chel}}}] \label{thm:CP}
    Let $X$ be a quasismooth Fano threefold weighted hypersurface with $\iota_X =1$. Then $X$ is birationally rigid.
\end{theorem}

\begin{Rem}
The following observations are made in \cite{CPR}: let $X$ be a quasismooth Fano threefold weighted hypersurface $X$ with $\iota_X = 1$. Suppose that $X$ is \textit{general} in its family. Then the following assertions hold.
\begin{enumerate}
\item any Sarkisov self-link $\sigma \colon X\dashrightarrow X$ of $X$ (if there are any) is of the form
\[
\xymatrix{
& \ar[ld]_{\varphi} Y \ar@{-->}[r]^{\tau} & Y \ar[rd]^{\varphi} & \\
X & & & X}
\]
where $\varphi$ is a divisorial extraction and $\tau$ is a flop, and in particular it does not involve an antiflip;
\item for any divisorial extraction $\varphi \colon Y \to X$, $\varphi$ is a maximal extraction if and only if $(-K_Y)^2 \in \operatorname{Int} \overline{\mathrm{NM}}_1 (Y)$ (\cite[Theorem~7.7.1]{CPR}).
Here, $\overline{\mathrm{NM}}_1 (Y) \subset \mathrm{N}_1 (Y)_{\mathbb{R}}$ is the \textbf{cone of nef curves} on $Y$ which is defined to be the dual of the cone $\overline{\mathrm{Eff}} (Y)$ of pseudo-effective divisors on $Y$. 
\end{enumerate}
By the work \cite{chel}, it turns out that both of the above observations hold true for \textit{every} quesismooth member of \textit{most} of the families, but do not hold true for some \textit{special} members of a few families. 

I explain this phenomena by picking family \textnumero~23 as an example.
Let $X = X_{14} \subset \mathbb{P} (1,2,3,4,5)$ be a quasismooth hypersurface.
Let $\mathbf{p} = \mathbf{p}_{x_2}$ be the singular point of type $\frac{1}{3} (1,1,2)$.
If $x_2^3 x_4$ appears in the defining equation with nonzero coefficient, then the point $\mathbf{p}$ is excluded as a maximal center (see \cite[Theorem~5.4.8]{CPR}).
The appearance of $x_2^3 x_4$ is one of the generality conditions assumed in \cite{CPR}.
Suppose that $x_2^3 x_4$ does not appear in the defining equation of $X$.
The point $\mathbf{p}$ is still not a maximal center if $x_2^2 x_3^2$ appears in the defining equation (see \cite[Page 68]{chel}).
Suppose that neither $x_2^3 x_4$ nor $x_2^2 x_3^2$ appear in the defining equation.
In this case, the point $\mathbf{p}$ is a maximal center and there is a Sarkisov self-link $\sigma \colon X \dashrightarrow X$ which sits in the diagram of the form
\[
\xymatrix{
& \ar[ld]_{\varphi} Y \ar@{-->}[r]^{\tau} & Y \ar[rd]^{\varphi} & \\
X \ar@{-->}[rrr]^{\sigma} & & & X}
\]
where $\varphi$ is the Kawamata blowup of $X$ at $\mathbf{p}$ and $\tau$ is the composite of antiflips, a flop, and flips (see \cite[Theorem~4.3.1]{chel}).
We explain the above diagram in more details.
Let $\Gamma$ and $\Delta$ be the proper transforms of the curves $(x_0 = x_1 = x_4 = 0) \subset X$ and $(x_0 = x_1 = x_3 = 0) \subset X$ on $Y$, respectively.
The proper transforms of the divisors $(x_0 = 0) \cap X$ and $(x_1 = 0) \cap X$ on $Y$ are linearly equivalent to $-K_Y$ and $-2K_Y$, respectively, and their intersection is the set $\Gamma \cup \Delta$.
Moreover, it is straightforward to compute $-K_Y \cdot \Gamma < 0$ and $- K_Y \cdot \Delta < 0$.
From these we deduce that $\Gamma$ and $\Delta$ are the only irreducible curves that intersect $-K_Y$ negatively.
It then turns out that the first contraction of the $2$-ray game constructing the map $\tau$ is the antiflipping contraction which contracts one of $\Gamma$ and $\Delta$ (further computations show that the contracted curve is $\Gamma$).
It follows that $X$ admits a Sarkisov self-link involving an antiflip, hence the observation (1) does not hold for this special member $X$.
Let $D \in \Cl (Y)$ be the proper transform of the Weil divisor $(x_0 = 0) \cap X \sim -K_X$ on $Y$.
Clearly $D$ is effective and we have $D \sim -K_Y$.
It follows that 
\[
D \cdot (-K_Y)^2 = (-K_Y)^3 = - \frac{1}{20} < 0.
\]
This shows that $(-K_Y)^2 \notin \overline{\mathrm{NM}}_1 (Y)$ and the observation (2) does not hold for $X$.
\end{Rem}
    
We next consider quasismooth Fano threefold WCIs in $\mathbb{P} (a_0, \dots, a_5)$ with $\iota_X = 1$.
There are exactly 85 such families and we index them by $I = \{1,\dots,85\}$ (see \cite[16.7]{CPR} and \cite{CCC}).
Among the $85$ families, there is a unique family consisting of smooth members.
It is the family of complete intersections of a quadric and a cubic in $\mathbb{P}^5$ and birational rigidity of its general members is proved in \cite{iskpuk}.
Corti and Mella \cite{cortimella} found first examples of Fano threefolds which are not birationally rigid but solid.
They proved that a general complete intersection of a cubic and a quartic in $\mathbb{P} (1,1,1,1,2,2)$ is birational to a factorial quartic threefold admitting a terminal singular point of type $cA_2$ and it is not birational to any other Mori fibre space.
Grinenko studied the family of the family of complete intersections of two cubics in $\mathbb{P} (1,1,1,1,1,2)$, and showed that any member in the family admits a Sarkisov link to a del Pezzo fibration over $\mathbb{P}^1$, which implies that it is not birationally solid.
\[
\begin{split}
I_{br} &= \{1, 8, 14, 20, 24, 31, 37, 45, 47, 51, 59, 60, 64, 71, 75, 76, 78, 84, 85\}, \\
I_{dP} &= \{2, 4, 5, 11, 12, 13\}, \\
I_F &= I \setminus I_{br} \cup I_{dP}.
\end{split}
\]
Based on these results, Okada started systematic studies for the remaining 82 families and obtained the following results:

\begin{Thm}[{\cite[Theorems~1.2 and 1.4]{okadaI}, \cite[Theorem~1.2]{okadaI}, \cite[Theorem~1.2]{okadaII}, \cite{okadaIII}, \cite[Theorem~1.1]{hamidmds}}] \label{thm:WCIcodim2index1}
Suppose $X$ is a quasismooth Fano threefold WCI of codimension $2$ and of index $1$ belonging to deformation family indexed by $I$.
\begin{enumerate}
\item Suppose that $X \in I_{br}$ and suppose in addition that $X$ is general if $X$ is a complete intersection of a quadric and cubic in $\mathbb{P}^5$.
Then $X$ is birationally rigid.
\item Suppose $X \in I_F$.
Then there exists a Sarkisov link $X \dashrightarrow X'$ to a $\mathbb{Q}$-factorial non-quasismooth Fano threefold weighted hypersurface $X'$ of index $1$.
Moreover, for a general $X$ belonging to suitable $35$ families out of $I_F$, any Mori fibre space which is birational to $X$ is either $X$ or $X'$, which implies $X$ is birationally solid but not birationally rigid.
\item If $X \in I_{dP}$, then there exists a Sarkisov link $X \dashrightarrow X'$ to a del Pezzo fibration $X'/\mathbb{P}^1$ over $\mathbb{P}^1$.
In particular, $X$ is not birationally solid.
\end{enumerate}
\end{Thm}

Thus, the study for quasismooth Fano threefold WCIs of codimension $2$ and of index $1$ has not been finished, and we need to consider the following conjecture and question.

\begin{Conj} \label{conj:codim2index1}
\begin{enumerate}
\item Any smooth commplete intersection of a quadric and a cubic in $\mathbb{P}^5$ is birationally rigid.
\item Any quasismooth member $X$ belonging to a family indexed by $I_F$ is not birational to a Mori fiber space other than $X$ and $X'$, where $X'$ is the Fano threefold given in (2) of Theorem~\ref{thm:WCIcodim2index1}.
\end{enumerate}
\end{Conj}

\begin{Quest}
Let $X$ be a quasismooth member of a family indexed by $I_{dP}$ and let $X'/\mathbb{P}^1$ is the del Pezzo fibration given in (3) of Theorem~\ref{thm:WCIcodim2index1}.
Then, is there a Mori fibre space other than $X$ and $X'/\mathbb{P}^1$ which is birational to $X$?
\end{Quest}

Finally, the only Fano threefold WCIs of codimension $\ge 3$ are the complete intersections of three quadrics in $\mathbb{P}^6$.
Let $X$ be such a Fano threefold.
Then it is known that $X$ contains a line and the blowup of $X$ along a line initiates a Sarkisov link to a conic bundle over $\mathbb{P}^2$ (see \cite[Theorem~4.3.3]{IskPro}), which implies that $X$ is not birationally solid.

\subsection{Fano index bigger than 1}
\label{sec:Fanoindex>1}

The case of higher Fano index is fundamentally different than the case of Fano index $1$. The following result allows to rule out birational solidity for many Fano threefold weighted complete intersections:
\begin{Thm}[{\cite[Corollary~2.3]{liviatiago}}] \label{thm:liviatiago}
Let $X$ be a  $\mathbb{Q}$-factorial terminal Fano $d$-fold with $d\geq 3$ and $A \in \Cl(X)$ a generator of the Class group of $X$. Consider the embedding given by the ring of sections of $A$
\begin{equation*}
    X \hookrightarrow \mathbb{P}(a_0,\ldots , a_N)
\end{equation*}
 Suppose that there are $a_i,\, a_j$ such that  $\lcm(a_i,a_j) < \iota_X$.\ Then, $X$ is not birationally solid.\ 
\end{Thm}

In \cite{hamidhyp}, the authors analyse the case of quasismooth hypersurfaces in weighted projective space with Fano index at least $2$. There are exactly 35 such families. They show that none of them is birationally rigid and, moreover, that only the following 5 families are possibly birationally solid:

\begin{multline} \label{eq:solidhyp}
    \Big\{X_{18}\subset \mathbb P(1,2,3,5,9),\, X_{22}\subset \mathbb P(1,2,3,7,11),\\ X_{26}\subset \mathbb P(1,2,5,7,13),\, 
    X_{38}\subset \mathbb P(2,3,5,11,19),\, X_{21}\subset \mathbb P(1,3,5,7,8)\Big\}
\end{multline}

They proceed to conjecture that these are indeed birationally solid which is confirmed in \cite{okadasolidhyp}. The following theorem is a combination of these two results:

\begin{Thm}[{\cite[Theorem~1.2]{hamidhyp},\cite[Theorem~1.3]{okadasolidhyp}}] \label{thm:highcodhyp}
Let $X$ be a quasismooth Fano threefold weighted hypersurface with $\iota_X \geq 2$. Then $X$ is not birationally rigid. Moreover, $X$ is birationally solid if and only if it belongs to one of the five families listed in Table \ref{eq:solidhyp}.
\end{Thm}

We add some more details of the results of \cite[Theorem~1.3]{okadasolidhyp} on the 5 families in \eqref{eq:solidhyp}.
Let $X$ be a quasismooth member of one of the $5$ families and let $\mathbf{p}$ be the (terminal quotient) singular point of the highest index.
Then it is proved in \cite{hamidhyp} that the Kawamata blowup $\varphi \colon Y \to X$ initiates a Sarkisov link to a Fano threefold hypersurface $\hat{X}$ of index $1$ in a weighted projective $3$-space (the construction of the links follows the method explained in \S~\ref{sec:SarkisovLinks}):
\[
\xymatrix{
& \ar[ld]_{\varphi} Y \ar@{-->}[r]^{\theta} & \hat{Y} \ar[rd]^{\hat{\varphi}} & \\
X & & & \hat{X}}
\]
The map $\theta$ is a birational map which is an isomorphism in codimension $1$.
The weighted hypersurface $\hat{X}$ is not quasismooth.
To be more precise, $\hat{X}$ admits a unique point $\hat{\mathbf{p}}$ at which $\hat{X}$ is not quasismooth.
The morphism $\hat{\varphi}$ is a divisorial contraction to the point $\hat{\mathbf{p}}$.
Let $X$ be a quasismooth member of family $110$.
Then, another model of $X$ is constructed in \cite{hamidhyp}.
Let $\mathbf{q} \in X$ be the singular point of type $\frac{1}{5}(1,1,4)$.
Then the Kawamata blowup of $X$ at $\mathbf{q}$ initiates a Sarkisov link $X \dashrightarrow \breve{X}$ to a Fano threefold weighted complete intersection $\breve{X} = \breve{X}_{6,7}$ of multidegree $(6,7)$ in $\mathbb{P} (1,1,2,2,3,5)$.
The WCI $\breve{X}$ has a unique point $\breve{\mathbf{q}}$ at which $\breve{X}$ is not quasismooth and the Sarkisov link $X \dashrightarrow \breve{X}$ terminates with the divisorial contraction to the point $\breve{\mathbf{q}} \in \breve{X}$.
The degree and weights of $\hat{X}$ and $\breve{X}_{6,7}$ as well as the types of singularities $\mathbf{p}, \hat{\mathbf{p}}, \mathbf{q}$ and $\breve{\mathbf{q}}$ are given in Table~\ref{table:solidmodel}. 
The following is the precise statement of \cite[Theorem~1.3]{okadasolidhyp}.
Note that the varieties $\hat{X}$ and $\breve{X}$ are all Mori fiber spaces, that is, they are all Fano threefolds with Picard number $1$.

\begin{itemize}
\item Let $X$ be a quasismooth member of families \textnumero $100, 101, 102$ and $103$.
Then the Fano threefolds $X$ and $\hat{X}$ are the only Mori fiber spaces in the birational equivalence class of $X$ (up to isomorphisms).
\item Let $X$ be a quasismooth member of families \textnumero $110$.
Then the Fano threefolds $X$, $\hat{X}$ and $\breve{X}$ are the only Mori fiber spaces in the birational equivalence class of $X$ (up to isomorphisms).
\end{itemize}

\begin{table}[t]
\centering
\caption{Birational models for birationally solid 5 families}
\label{table:solidmodel} 
\renewcommand{\arraystretch}{1.2}
\begin{tabular}[]{ccccc}
\toprule
\textnumero & $X_d \subset \mathbb{P} (a_0,\dots,a_4)$ & $\mathbf{p}$ (or $\mathbf{q}$) & $\hat{X}_{e} \subset \mathbb{P} (b_0, b_1, b_2, b_3, b_4)$ & $\hat{\mathbf{p}}$ (or $\hat{\mathbf{q}}$) \\
\midrule
100 & $X_{18} \subset \mathbb{P} (1,2,3,5,9)$ & $\frac{1}{5} (1,2,3)$ & $\hat{X}_{10} \subset \mathbb{P} (1,1,1,3,5)$ & $cE_6$ \\
101 & $X_{22} \subset \mathbb{P} (1,2,3,7,11)$ & $\frac{1}{7} (1,2,5)$ & $X_{12} \subset \mathbb{P} (1,1,1,4,6)$ & $cE_{7,8}$ \\
102 & $X_{26} \subset \mathbb{P} (1,2,5,7,13)$ & $\frac{1}{7} (1,3,4)$ & $\hat{X}_{14} \subset \mathbb{P} (1,1,2,4,7)$ & $cE/2$ \\
103 & $X_{38} \subset \mathbb{P} (2,3,5,11,19)$ & $\frac{1}{11}(1,4,7)$ & $\hat{X}_{22} \subset \mathbb{P}(1,1,3,7,11)$ & $cE_8$ \\
110 & $X_{21} \subset \mathbb{P} (1,3,5,7,8)$ & $\frac{1}{8} (1,3,5)$ & $\hat{X}_7 \subset \mathbb{P}(1,1,1,2,3)$ & $cE_7$ \\
& & $\frac{1}{5}(1,1,4)$ & $\breve{X}_{6,7} \subset \mathbb{P}(1,1,2,2,3,5)$ & $cD/2$ \\
\bottomrule
\end{tabular}
\end{table}

 A similar situation holds for $X$ a Fano threefold embedded in a weighted projective space as a codimension $2$ complete intersection with index larger than 1. There are exactly 40 such families and we index them by $I=\{86, \ldots,125\}$. See \cite[The Big Table]{tiagononrigid} for the definition of each family in terms of $I$. Moreover, let
\begin{align*}
I_{Cb}&=\{87,112,113,118,119\}\\
I_{dP}&=\{88,89,90,91,103,114,116,120,121,122,123,124,125\}
\end{align*}
and
\begin{align*}
I_{nS}=I_{Cb} \cup I_{dP}, \qquad
I_{S}=I \setminus I_{nS}.
\end{align*}

The following is the main Theorem of \cite{tiagononrigid}:

\begin{Thm}[{\cite[Theorem~1.3]{tiagononrigid}}] \label{thm:maintiago}
Suppose $X$ is a quasismooth Fano threefold WCI of codimension $2$ and of index $\geq 2$ belonging to deformation family indexed by $I$. 
\begin{enumerate}
	\item If $X \in I_{nS}$ then $X$ is not birationally solid. In particular,
	\begin{enumerate}
		\item if $X \in I_{Cb}$ then there is a Sarkisov link from $X$ to a conic bundle $Y/S$, where $S$ is smooth or has $A_1$ or $A_2$ singularities.
		\item if $X \in I_{dP}$ then there is a Sarkisov link from $X$ to a del Pezzo fibration whose generic fibre has degree 1, 2 or 3.
	\end{enumerate}
		\item If $X \in I_{S}$ then there is a Sarkisov link from $X$  to a terminal non-quasismooth weighted hypersurface or codimension 2 Fano threefold.
\end{enumerate}
In particular, $X$ is not birationally rigid. 
\end{Thm}

All the information regarding the decomposition of all the Sarkisov links from $X$ as above can be seen in \cite[Section: The Big Table]{tiagononrigid}. On the other hand, the picture regarding birational solidity of these families is wide open. The following is a natural question:

\begin{Quest}[Duarte Guerreiro] \label{conj:solid}
Is a general member $X$ in $I_S$ birationally solid ?
\end{Quest}

The affirmative answer to Question~\ref{conj:solid} is known only to 1 family: it is proved in \cite{Okasolidcod2ind2} that  a quasismooth $X_{12,14} \subset \mathbb{P} (1,2,3,4,7,11)$, that is, a quasismooth member of family \textnumero~100, is birationally solid.

\section{Some remarks on non-complete intersection Fano threefolds}
\label{sec:nonWCI}

\subsection{Codimension 3} Let $X$ be a quasismooth Fano threefold embedded in a weighted projective space as a codimension $3$ subvariety such that $\iota_X = 1$. Then either $X$ is a complete intersection of $3$ quadrics in $\mathbb P^6$ or $X$ is a \emph{Pfaffian}, see \cite{pfaffians}. Such Fano threefolds are given by the $4\times 4$ Pfaffians of a skew $5\times 5$ matrix by the structure theorems of Buchsbaum and
Eisenbud, see \cite{buchsbaum1977algebra}. In total there are 70 deformation families. Birational rigidity in the case of Pfaffian Fano threefolds of index $1$ has been characterized in \cite{hamidokadapfaff}.  If the Fano index of $X$ is at least $2$, then $X$ is not birationally solid by Theorem \ref{thm:liviatiago}.

\subsection{Codimension 4}  Let $X$ be a quasismooth Fano threefold embedded with codimension $4$ in a weighted projective space. Then there are no structure theorems for $X$ as for the lower codimensions. However a first strategy was introduced in \cite{cod4tomJerry}, see also \cite{cod4tomJerryII}, based on the Kustin-Miller unprojection which allowed the construction of more than a hundred families of anticanonically polarised Fano threefolds embedded in codimension $4$ in weighted projective space.  Recently, Coughlan and Ducat give a completely different construction of prime Fano threefolds embedded as codimension $4$ in a weighted projective arising from cluster algebras, see \cite{coughlan2020constructing}. In \cite{cod4Okada}, birational superrigidity of cluster Fano threefolds is characterised.
Takagi generalizes the method and constructs further Fano threefolds embedded as codimension 4 in a weighted projective space in a series of papers such as \cite{TakagiKey1}, \cite{TakagiKey2}, \cite{TakagiKey3} and \cite{TakagiKey4}.

For Fano index $2$ a construction has appeared in \cite{liviahalfelephants} and birational geometry has been analysed in \cite[Theorem~1.2]{liviatiago} where it is shown, in particular, that they are all birationally non-rigid.

\section{Some remarks on non-quasismooth Fano threefolds}
\label{sec:nonqsm}

It is expected that the presence of singularities on a given Fano variety $X$ will increase the number of Mori fibre spaces in the birational class of $X$. The first such example is the one from Corti and Mella:

\begin{Ex}[{\cite{cortimella}}]  \label{ex:cortimella}
The classical example is of a quartic threefold $X$ with a unique singular point analytically equivalent to $(xy+z^3+t^3=0) \subset \mathbb C^4$, which is a $cA_2$ singularity. Then $X$ is factorial and terminal \cite[Lemma~6.4.4]{kawakitabook}. Let $Y/S$ be a Mori fibre space. If $X$ and $Y$ are birational then $S$ is a point and either $Y$ is biregular to $X$ or $Y$ is biregular to $Y_{3,4} \subset \mathbb P(1,1,1,1,2,2)$, the general complete intersection of a cubic and a quartic. See Example \ref{ex:quarticsarkisovlink} for details on the construction of the Sarkisov link to $Y_{3,4}$. Moreover, blowing up a line through the singular point induces a birational involution of $X$ while a weighted $(2,1,1,1)$-blowup of the singular point gives a birational map to $Y_{3,4}$. 
\end{Ex}

  In \cite{hamidquartic} the authors confirm this behaviour by constructing factorial quartic threefolds which are not birationally rigid. They showed the following
    \begin{theorem}[{\cite[Theorem~1.2]{hamidquartic}}]
      Let $X$ be a factorial quartic threefold with only $cA_n$ singularities. Then $n\leq 7$. Moreover, for each such $n\geq 2$ there is a non-birationally rigid quartic threefold with a $cA_n$ singularity.
  \end{theorem}

  However, it is proved in \cite{PukhlikovquarticcA1} that a quartic threefold admitting a unique singular point which is of type $cA_1$ is birationally rigid.
  The following conjecture remains open:

\begin{Conj}
    A $\mathbb{Q}$-factorial (or equivalently factorial) quartic threefold with only $cA_1$ singular points is birationally rigid. 
\end{Conj}

Similarily, sextic double solids, i.e., double covers of $\mathbb P^3$ ramified along sextic surfaces are birationally superrigid whenever the sextic surface is smooth. See \cite{iskovskikh1979birational}. On the other hand, the introduction of singularities breaks rigidity. Indeed by \cite[Theorem~B]{Erik}, up to some explicit generality conditions, every factorial sextic double solid with a $cA_n$ singularity where $n\geq 4$ is not birationally rigid.   

\begin{Conj}[Paemurru]
    A $\mathbb{Q}$-factorial (or equivalently factorial) sextic double solid with only $cA_1$ and $cA_2$ singular points is birationally superrigid. A $\mathbb{Q}$-factorial sextic double solid with only $cA_1$, $cA_2$ and $cA_3$ singular points is birationally rigid.
\end{Conj}

The superrigidity of $\mathbb{Q}$-factorial sextic double solids with only $cA_1$ singular points has been confirmed in \cite[Theorem~4.1]{2n2-ineq}.

\section{Some more questions on birational rigidity and solidigy}
We have completed overall explanations of known results of birational geometry of Fano threefolds related to birational rigidity in \S~\ref{sec:bigpicture}, \S~\ref{sec:nonWCI} and \S~\ref{sec:nonqsm}.
We give some questions regarding birational rigidity and solidity.

From the known results, we can vaguely say that Fano index, codimension and singularities can be thought of as obstructions to birational rigidity/solidity.

\begin{Quest}
   Find natural numbers $i_{\mathrm{br}}$, $i_{\mathrm{bs}}$, $c_{\mathrm{br}}$ and $c_{\mathrm{bs}}$ such that for each Fano threefold $X$ of Picard rank $1$, we have the following implications
   \begin{itemize}
       \item $\iota_X \geq i_{\mathrm{br}}$ \quad (resp.  $\mathrm{cod}_X \geq c_{\mathrm{br}}$) \quad $\implies$ \quad $X$ is not birationally rigid.
       \item $\iota_X \geq i_{\mathrm{bs}}$ \quad (resp.  $\mathrm{cod}_X \geq c_{\mathrm{bs}}$) \quad $\implies$ \quad $X$ is not birationally solid.
   \end{itemize}
\end{Quest}

Note that it is proved by Prokhorov \cite{ProkhorovRatLargeindex} that a Fano threefold $X$ is rational if $\iota_X \geq 8$ in which the lower bound $8$ is optimal.
Currently we only know that $i_{\mathrm{br}} \geq 2$, $i_{\mathrm{bs}} \geq 4$, $c_{\mathrm{br}} \geq 5$ and $c_{\mathrm{bs}} \geq 5$.

\begin{Quest}
    Is there a birationally rigid (resp.\ solid) Fano threefold $X$ with $\iota_X \geq 2$ (resp.\ $\iota_X \geq 4$)?
    Is there a birationally rigid or solid Fano threefold $X$ with $\mathrm{cod}_X \geq 5$?
\end{Quest}

For a birationally solid Fano variety $X$, the set
\[
\mathcal{P} (X) := \{X' \mid \text{$X'$ is a Fano variety of Picard rank $1$ such that $X' \sim_{\mathrm{bir}} X$}\}/\cong
\]
is called the \textbf{pliability set} of $X$.
In all the known cases, any model $X' \in \mathcal{P} (X) \setminus \{X\}$ of a given birationally solid but not birationally rigid Fano threefold $X$ has a non-quotient singularity, that is, $X'$ is not quasismooth even if $X$ is quasismooth.
Moreover, there is always a Fano variety $X' \in \mathcal{P} (X)$ whose Fano index is $1$. 

\begin{Quest} Suppose $X$ is birationally solid but not birationally rigid. Is it possible that all Fanos in the set $\mathcal{P} (X)$ are quasismooth ?
Is it possible that all Fanos in $\mathcal{P} (X)$ have Fano index at least $2$?
\end{Quest}

For any known example of a birationally solid Fano threefold $X$, we have $|\mathcal{P} (X)| \leq 3$ and, moreover, any model $X' \in \mathcal{P} (X) \setminus \{X'\}$ is connected to $X$ by a single Sarkisov link, that is, the \textbf{pliability graph} of $X$, which is the graph whose vertices are the elements of $\mathcal{P} (X)$ and two distinct vertices $X_1, X_2 \in \mathcal{P} (X)$ are connected by an edge if and only if there is a Sarkisov link between $X_1$ and $X_2$, is complete. 

\begin{Quest}
\label{quest:Pli}
\begin{enumerate}
\item For an integer $n \geq 3$, does there exist a positive integer $N (n)$ depending only on $n$ such that $\mathcal{P} (X) \leq N (n)$ for any birationally solid Fano variety $X$ of dimension $n$?
\item Does there exist a birationally solid Fano variety $X$ such that $\mathcal{P} (X)$ is an infinite set?
\end{enumerate}
\end{Quest}

Very recently, the existence of a del Pezzo surface $S$ of Picard number $1$ with infinite pliability set $\mathcal{P} (S)$ is proved in \cite[Theorem~C]{KurzEgor}.
The surface $S$ is defined over an algebraically non-closed field, hence Question~\ref{quest:Pli} remains valid over an algebraically closed base field.

\begin{Quest}
Is there a birationally solid Fano threefold $X$ whose pliability graph is not complete?
\end{Quest}

We can define the pliability set $\mathcal{P} (X/S)$ and its associated graph for an arbitrary Mori fiber space $X/S$ in a similar way, see for example \cite[Definition 1.5]{cortimella} for the definition. 
There is a Fano threefold $X$ of Picard number $1$ for which $|\mathrm{Pli} (X)| = 9$ and its pliability graph is star-shaped and is not complete \cite[\S 3.3]{Sarikyan}: more precisely, there are del Pezzo fibrations $X_i/\mathbb{P}^1$ for $1 \leq i \leq 8$ such that each $X_i/\mathbb{P}^1$ is connected to $X$ by a Sarkisov link while $X_i/\mathbb{P}^1$ cannot be connected to $X_j/\mathbb{P}^1$ by a Sarkisov link for any $i \ne j$.
Note that this Fano threefold $X$ is not birationally solid.

\section{Rationality problems of Fano threefold WCIs}

We discuss the rationality problem of quasismooth Fano threefold WCIs. 

Any smooth complete intersection of three quadrics in $\mathbb{P}^6$ is proved to be irrational by analyzing their intermediate Jacobians \cite{3quad}.
To the best knowledge of the authors, the only results on rationality questions for Fano threefold WCIs of codimension $2$ are due to birational rigidity and solidity except for smooth complete intersections of two quadrics in $\mathbb{P}^5$ which are known to be rational, hence we have nothing to add beyond the results in \S~\ref{sec:bigpicture}.
We concentrate on the rationality problem of quasismooth Fano threefold weighted hypersurfaces.

By Theorem~\ref{thm:CP}, any quasismooth Fano threefold weighted hypersurface of index $1$ is irrational since it is birationally rigid.
It remains to consider quasismooth Fano threefold weighted hypersurfaces of index $> 1$ which form $35$ deformation families.
We have the following criterion for rationality (see \cite{EsserRational} for other rationality criteria).
We say that a weighted hypersurface $X$ in a weighted projective space $\mathbb{P} (a_0, \dots, a_{n+1})$ with homogeneous coordinates $x_0, \dots, x_{n+1}$ is a \textbf{weighted cone} if there is $i \in \{0, \dots, n+1\}$ such that the defining equation does not involve the variable $x_i$.

\begin{Lem}[{\cite[Lemma~6.1]{okadastratFano}}] \label{lem:ratcriWH}
Let $a_0, a_1, \dots, a_{n+1}$ and $d$ be positive integers such that $a_0 \leq \cdots \leq a_{n+1}$.
Suppose that either $d < 2 a_{n+1}$ or $d = 2 a_{n+1} = 2 a_n$.
Then a weighted hypersurface of degree $d$ in $\mathbb{P} (a_0, \dots, a_{n+1})$ is rational if it is irreducible, reduced, and not a weighted cone.
\end{Lem}

\begin{table}[t]
\centering
\caption{Rationality of quasismooth Fano threefold weighted hypersurfaces of index $\geq 2$: 
In the column ``Rat'', the sign $-$ means that any quasismooth member of the family is irrational.
The column ``Ref'' indicates the reference.}
\label{table:ratWH} 
\begin{tabular}[]{cccccccc}
\toprule
\textnumero & $X_d \subset \mathbb{P} (a_0,\dots,a_4)$ & Rat & Ref & \textnumero & $X_d \subset \mathbb{P} (a_0,\dots,a_4)$ & Rat & Ref \\
\midrule
96 & $X_3 \subset \mathbb{P} (1,1,1,1,1)$ & $-$ & \cite{CG} & 107 & $X_6 \subset \mathbb{P} (1, 1, 2, 2, 3)$ & $-$ & \cite{ProkhorovWH}  \\
97 & $X_4 \subset \mathbb{P} (1,1,1,1,2)$ & $-$ & \cite{voisinquarticds} & 108 & $X_{12} \subset \mathbb{P} (1, 2, 3, 4, 5)$ & ? \\
98 & $X_6 \subset \mathbb{P} (1,1,1,2,3)$ & $-$ & \cite{GrinenkoVeronese} & 109 & $X_{15} \subset \mathbb{P} (1, 2, 3, 5, 7)$ & ? \\
99 & $X_{10} \subset \mathbb{P} (1,1,2,3,5)$ & ? & & 110 & $X_{21} \subset \mathbb{P} (1,3,5,7,8)$ & $-$ & \cite{okadasolidhyp} \\
100 & $X_{18} \subset \mathbb{P} (1,2,3,5,9)$ & $-$ & \cite{okadasolidhyp} & 116 & $X_{10} \subset \mathbb{P} (1,2,3,4,5)$ & $-$ & \cite{ProkhorovbirgeomII} \\
101 & $X_{22} \subset \mathbb{P} (1,2,3,7,11)$ & $-$ & \cite{okadasolidhyp} & 117 & $X_{15} \subset \mathbb{P} (1, 3, 4, 5, 7)$ & ? & \\
102 & $X_{26} \subset \mathbb{P} (1,2,5,7,13)$ & $-$ & \cite{okadasolidhyp} & 122 & $X_{14} \subset \mathbb{P} (2,3,4,5,7)$ & ? & \\
103 & $X_{38} \subset \mathbb{P} (2,3,5,11,19)$ & $-$ & \cite{okadasolidhyp} & \\
\bottomrule
\end{tabular}
\end{table}

There are $20$ families out of $35$ which satisfy the assumptions of Lemma~\ref{lem:ratcriWH}, hence any quasismooth member in each of the $20$ families is rational.
In Table~\ref{table:ratWH}, we summarize the known results on the rationality questions for the remaining $15$ families.
As discussed in \S~\ref{sec:Fanoindex>1}, any quasismooth member in families \textnumero~100, 101, 102, 103 and 110 are birationally solid (see Theorem~\ref{thm:highcodhyp}), hence it is irrational.
Any quasismooth member in the remaining 10 families thus admits a strict Mori fiber space, which is either a del Pezzo fibration or a conic bundle, as its birational model.
Prokhorov studied in detail a conic bundle structure of a quasismooth member of families 107 and 116 in \cite{ProkhorovWH} and \cite{ProkhorovbirgeomII}, respectively, and showed that the associated Prym variety cannot be isomorphic as a principally polarized abelian variety to a sum of Jacobian of curves, which implies its irrationality.
To the best knowledge of the authors, the rationality questions for families \textnumero~99, 108, 109, 117 and 122 remain open.

\begin{Quest}
\label{quest:irratWH}
Are quasismooth members in families \textnumero 99, 108, 109, 117 and 122 irrational?
\end{Quest}

In the case of weighted hypersurfaces, the construction of rationality is given by Lemma~\ref{lem:ratcriWH} and this construction is given simply by observing that a projection from a singular point of the highest index is a birational map onto a weighted projective $3$-fold.
This simple construction of rationality is not applicable to quasismooth Fano threefold WCIs of codimension $2$ and of index $\geq 2$, and currently we do not know the existence of rational quasismooth Fano threefold WCIs of codimension $2$ other than complete intersections of two quadrics in $\mathbb{P}^5$.

\begin{Quest}
Are there quasismooth Fano threefold WCIs of codimension $2$ and of index $\geq 2$ which are rational other than complete intersections of two quadrics in $\mathbb{P}^5$?
\end{Quest}

\begin{Rem}
There is a weaker notion of rationality: a variety $X$ is said to be \textbf{stably rational} if there is a nonnegative integer $m$ such that $X \times \mathbb{P}^m$ is rational, where $\mathbb{P}^0$ stands for $\operatorname{Spec} \mathbb{C}$.
We have the implications:
\[
\text{rational} \ \Longrightarrow \ \text{stably rational} \ \Longrightarrow \ \text{unirational} \ \Longrightarrow \text{rationally connected}
\]
The above implications cannot be reversed except possibly for the last one.
See \cite{BCTSSD} for the existence of a stably rational variety which is not rational.
See \cite{Voisindegeneration} for the existence of unirational varieties which are not stably rational.
In \cite{okadastratFano}, stable rationality of quasismooth Fano threefold weighted hypersurfaces are studied by the specialization method which is originally invented by Voisin \cite{Voisindegeneration} and developed by Koll\'ar \cite{Kollarnonrathyp}, Totaro \cite{Totarostablyrational}, Colliot-Th\'el\`ene--Pirutka \cite{CTPirutka}.
As a consequence, it is proved that the following Fano threefolds are not stably rational:
\begin{itemize}
\item Very general Fano threefold weighted hypersurfaces of index $1$.
\item Very general members in each of the $15$ families listed in Table~\ref{table:ratWH}.
\end{itemize}
Based on these results, we believe the positive answer to Question~\ref{quest:irratWH}.
\end{Rem}

\section{Relation to K-stability}

The following beautiful conjecture can be found for example in \cite[Conjecture~2.1]{newdir}:
\begin{Conj}
    Let $X$ be a smooth uniruled variety. Then $X$ is birational to a Mori fibre space where the general fibre is K-stable.
\end{Conj}

In light of the above conjecture it is expected, in particular, that birationally rigid Fano varieties be K-stable. The following result is a first step in that direction.

\begin{Def}
    Let $X$ be a normal projective $\mathbb{Q}$-Gorenstein variety with ample anticanonical divisor. 
    The \textbf{alpha invariant} of $X$ denoted $\alpha(X)$ is defined as the supremum of all $t>0$ such that $(X,D)$ is log canonical for every effective $\mathbb Q$-divisor $D\sim_{\mathbb Q} -K_X$.
\end{Def}

\begin{Thm}[{\cite[Theorem~1.2]{Stibitz_Zhuang_2019}}] \label{thm:superrigalpha}
   Let $X$ be a Fano variety of Picard number $1$. If $X$ is birationally superrigid and $\alpha(X) \geq \frac{1}{2}$ (resp, $> \frac{1}{2}$), then $X$ is K-semistable (resp. K-stable).
\end{Thm}

\begin{Rem}
    In \cite{kim2023k}, the authors use Theorem \ref{thm:superrigalpha} to establish K-stability of \emph{birationally superrigid} quasismooth Fano threefold weighted hypersurfaces of index $1$.
\end{Rem}

Theorem \ref{thm:superrigalpha} can be slightly improved. Indeed if $X$ is superrigid and $\alpha(X)\geq 1/2$, then $X$ is uniformely K-stable. See \cite[Theorem~9.6]{xu2021k}. On the other hand, every known birational superrigid Fano variety satisfies $\alpha(X)>\frac{1}{2}$.
Note that also there are birationally superrigid Fano varieties whose alpha invariants are arbitrary close to $1/2$ \cite[Theorem~1.4]{CSZalpha}: more precisely, for any real number $\varepsilon > 0$, there is a birationally superrigid Fano variety $X$ (with only terminal singularities) for which $1/2 \leq \alpha (X) \leq 1/2 + \varepsilon$.

\begin{Quest}[{cf.\cite[Question~1.5]{Stibitz_Zhuang_2019}}]
Is it true that $\alpha (X) \geq 1/2$ for any birationally superrigid Fano variety $X$?
Is there an example which achieves the equality?
\end{Quest}

After the works \cite{kim2023k} and \cite{SanoTasinDelta}, the problem of K-stability for quasismooth Fano threefold weighted hypersurfaces of Fano index $1$ is settled in \cite{campookada}. Their main result is the following:

\begin{Thm}[{\cite[Theorem~1.5]{campookada}}]
    Any quasismooth Fano threefold weighted hypersurface of Fano index 1 is K-stable. 
\end{Thm}

We end this section with the following 
\begin{Quest} \label{quest:qsK-stab}
    Let $X$ be a quasismooth Fano threefold WCI of Fano index $1$. Is $X$ K-stable ?
\end{Quest}

\begin{Rem}
    The assumption on the Fano index in Question \ref{quest:qsK-stab} is necessary. For example, a smooth quadric hypersurface is strictly K-semistable. See for instance \cite[Corollary~4.11]{AZI} for an algebraic proof.
\end{Rem}

\section{Summary}

We present a summary of the state of the art on the problem of determining birational superrigidity, birational rigidity and birational solidity in the context of quasismooth Fano threefold complete intersections. 

\begin{itemize}
\item In the Tables~\ref{table:WCIcod1ind1}, \ref{table:WCIcod2ind1} and \ref{table:WCIcod1indhigh}, the column ``BP" stands for a birational property of quasismooth members of the corresponding family.
It indicates whether quasismooth members are birationally rigid ($\mathrm{BR}$), birationally solid ($\mathrm{BS}$), or it is birational to a del Pezzo fibration of degree $1$ over $\mathbb{P}^1$($\mathrm{dP}_1$).
\item In Table~\ref{table:WCIcod2ind1}, any quasismooth member of a family whose ``BP" column is blank is expected to be birationally solid but currently this question remains open.
Moreover, $\mathrm{BR}_*$ and $\mathrm{BS}_*$ mean that birational rigidity and solidity, respectively, of quasismooth members of the corresponding family are proven only under a generality assumption.  
\item In Table~\ref{table:WCIcod1indhigh}, the column ``BS" and ``Rat" stand for birational solidity and rationality, respectively, in which the sign $+/-$ means that any quasismooth member is birationally solid/not birationally solid and rational/irrational, respectively.
The sign ``?" means that the rationality (of any quasismooth member) is unknown. 
\item In the column ``Models" of Table~\ref{table:WCIcod2highind}, it is indicated a Mori fibre space that is connected to a member of the corresponding family by a Sarkisov link initiated by the Kawamata blowup at the singular point given in the column ``Centre".
Here, $Y/S$ and $\mathrm{dP}_i$ stand for a conic bundle over $S$ and del Pezzo fibration of degree $i$ over $\mathbb{P}^1$, respectively, and $X_{d_1, d_2} \subset \mathbb{P}(b_0,\dots,b_5)/\bm{\mu}_2$ stands for the quotient of a WCI $X_{d_1,d_2} \subset \mathbb{P}(b_0,\dots,b_5)$ by a suitable action of $\bm{\mu}_2$, see \cite[\S~6.2]{tiagononrigid}.
\end{itemize}

\begin{table}[t]
\centering
\caption{Fano threefold weighted hypersurfaces of index $1$}
\label{table:WCIcod1ind1} 
\begin{tabular}{cccccc}
\toprule
\textnumero & $X_d \subset \mathbb{P} (a_0,\dots,a_4)$ & BP & \textnumero & $X_d \subset \mathbb{P} (a_0,\dots,a_4)$ & BP \\
\midrule
1 & $X_4 \subset \mathbb{P} (1,1,1,1,1)$ & $\mathrm{BSR}$ & 49 & $X_{21} \subset \mathbb{P} (1,3,5,6,7)$ & $\mathrm{BSR}$ \\
2 & $X_5 \subset \mathbb{P} (1,1,1,1,2)$ & $\mathrm{BR}$ & 50 & $X_{22} \subset \mathbb{P} (1,1,3,7,11)$ & $\mathrm{BSR}$ \\
3 & $X_6 \subset \mathbb{P} (1,1,1,1,3)$ & $\mathrm{BSR}$ & 51 & $X_{22} \subset \mathbb{P} (1,1,4,6,11)$ & $\mathrm{BSR}$ \\
4 & $X_6 \subset \mathbb{P} (1,1,1,2,2)$ & $\mathrm{BR}$ & 52 & $X_{22} \subset \mathbb{P} (1,2,4,5,11)$ & $\mathrm{BSR}$ \\
5 & $X_7 \subset \mathbb{P} (1,1,1,2,3)$ & $\mathrm{BR}$ & 53 & $X_{24} \subset \mathbb{P} (1,1,3,8,12)$ & $\mathrm{BSR}$ \\
6 & $X_8 \subset \mathbb{P} (1,1,1,2,4)$ & $\mathrm{BR}$ & 54 & $X_{24} \subset \mathbb{P} (1,1,6,8,9)$ & $\mathrm{BR}$ \\
7 & $X_8 \subset \mathbb{P} (1,1,2,2,3)$ & $\mathrm{BR}$ & 55 & $X_{24} \subset \mathbb{P} (1,2,3,7,12)$ & $\mathrm{BSR}$ \\
8 & $X_9 \subset \mathbb{P} (1,1,1,3,4)$ & $\mathrm{BR}$ & 56 & $X_{24} \subset \mathbb{P} (1,2,3,8,11)$ & $\mathrm{BR}$ \\
9 & $X_9 \subset \mathbb{P} (1,1,2,3,3)$ & $\mathrm{BR}$ & 57 & $X_{24} \subset \mathbb{P} (1,3,4,5,12)$ & $\mathrm{BSR}$ \\
10 & $X_{10} \subset \mathbb{P} (1,1,1,3,5)$ & $\mathrm{BSR}$ & 58 & $X_{24} \subset \mathbb{P} (1,3,4,7,10)$ & $\mathrm{BR}$ \\
11 & $X_{10} \subset \mathbb{P} (1,1,2,2,5)$ & $\mathrm{BSR}$ & 59 & $X_{24} \subset \mathbb{P} (1,3,6,7,8)$ & $\mathrm{BSR}$ \\
12 & $X_{10} \subset \mathbb{P} (1,1,2,3,4)$ & $\mathrm{BR}$ & 60 & $X_{24} \subset \mathbb{P} (1,4,5,6,9)$ & $\mathrm{BR}$ \\
13 & $X_{11} \subset \mathbb{P} (1,1,2,3,5)$ & $\mathrm{BR}$ & 61 & $X_{25} \subset \mathbb{P} (1,4,5,7,9)$ & $\mathrm{BR}$ \\
14 & $X_{12} \subset \mathbb{P} (1,1,1,4,6)$ & $\mathrm{BSR}$ & 62 & $X_{26} \subset \mathbb{P} (1,1,5,7,13)$ & $\mathrm{BSR}$ \\
15 & $X_{12} \subset \mathbb{P} (1,1,2,3,6)$ & $\mathrm{BR}$ & 63 & $X_{26} \subset \mathbb{P} (1,2,3,8,13)$ & $\mathrm{BSR}$ \\
16 & $X_{12} \subset \mathbb{P} (1,1,2,4,5)$ & $\mathrm{BR}$ & 64 & $X_{26} \subset \mathbb{P} (1,2,5,6,13)$ & $\mathrm{BSR}$ \\
17 & $X_{12} \subset \mathbb{P} (1,1,3,4,4)$ & $\mathrm{BR}$ & 65 & $X_{27} \subset \mathbb{P} (1,2,5,9,11)$ & $\mathrm{BR}$ \\
18 & $X_{12} \subset \mathbb{P} (1,2,2,3,5)$ & $\mathrm{BR}$ & 66 & $X_{27} \subset \mathbb{P} (1,5,6,7,9)$ & $\mathrm{BSR}$ \\
19 & $X_{12} \subset \mathbb{P} (1,2,3,3,4)$ & $\mathrm{BSR}$ & 67 & $X_{28} \subset \mathbb{P} (1,1,4,9,14)$ & $\mathrm{BSR}$ \\
20 & $X_{13} \subset \mathbb{P} (1,1,3,4,5)$ & $\mathrm{BR}$ & 68 & $X_{28} \subset \mathbb{P} (1,3,4,7,14)$ & $\mathrm{BR}$ \\
21 & $X_{14} \subset \mathbb{P} (1,1,2,4,7)$ & $\mathrm{BSR}$ & 69 & $X_{28} \subset \mathbb{P} (1,4,6,7,11)$ & $\mathrm{BR}$ \\
22 & $X_{14} \subset \mathbb{P} (1,2,2,3,7)$ & $\mathrm{BSR}$ & 70 & $X_{30} \subset \mathbb{P} (1,1,4,10,15)$ & $\mathrm{BSR}$ \\
23 & $X_{14} \subset \mathbb{P} (1,2,3,4,5)$ & $\mathrm{BR}$ & 71 & $X_{30} \subset \mathbb{P} (1,1,6,8,15)$ & $\mathrm{BSR}$ \\
24 & $X_{15} \subset \mathbb{P} (1,1,2,5,7)$ & $\mathrm{BR}$ & 72 & $X_{30} \subset \mathbb{P} (1,2,3,10,15)$ & $\mathrm{BSR}$ \\
25 & $X_{15} \subset \mathbb{P} (1,1,3,4,7)$ & $\mathrm{BR}$ & 73 & $X_{30} \subset \mathbb{P} (1,2,6,7,15)$ & $\mathrm{BSR}$ \\
26 & $X_{15} \subset \mathbb{P} (1,1,3,5,6)$ & $\mathrm{BR}$ & 74 & $X_{30} \subset \mathbb{P} (1,3,4,10,13)$ & $\mathrm{BR}$ \\
27 & $X_{15} \subset \mathbb{P} (1,2,3,5,5)$ & $\mathrm{BR}$ & 75 & $X_{30} \subset \mathbb{P} (1,4,5,6,15)$ & $\mathrm{BSR}$ \\
28 & $X_{15} \subset \mathbb{P} (1,3,3,4,5)$ & $\mathrm{BSR}$ & 76 & $X_{30} \subset \mathbb{P} (1,5,6,8,11)$ & $\mathrm{BR}$ \\
29 & $X_{16} \subset \mathbb{P} (1,1,2,5,8)$ & $\mathrm{BSR}$ & 77 & $X_{32} \subset \mathbb{P} (1,2,5,9,16)$ & $\mathrm{BSR}$ \\
30 & $X_{16} \subset \mathbb{P} (1,1,3,4,8)$ & $\mathrm{BR}$ & 78 & $X_{32} \subset \mathbb{P} (1,4,5,7,16)$ & $\mathrm{BSR}$ \\
31 & $X_{16} \subset \mathbb{P} (1,1,4,5,6)$ & $\mathrm{BR}$ & 79 & $X_{33} \subset \mathbb{P} (1,3,5,11,14)$ & $\mathrm{BR}$ \\
32 & $X_{16} \subset \mathbb{P} (1,2,3,4,7)$ & $\mathrm{BR}$ & 80 & $X_{34} \subset \mathbb{P} (1,3,4,10,17)$ & $\mathrm{BSR}$ \\
33 & $X_{17} \subset \mathbb{P} (1,2,3,5,7)$ & $\mathrm{BR}$ & 81 & $X_{34} \subset \mathbb{P} (1,4,6,7,17)$ & $\mathrm{BSR}$ \\
34 & $X_{18} \subset \mathbb{P} (1,1,2,6,9)$ & $\mathrm{BSR}$ & 82 & $X_{36} \subset \mathbb{P} (1,1,5,12,18)$ & $\mathrm{BSR}$ \\
35 & $X_{18} \subset \mathbb{P} (1,1,3,5,9)$ & $\mathrm{BSR}$ & 83 & $X_{36} \subset \mathbb{P} (1,3,4,11,18)$ & $\mathrm{BSR}$ \\
36 & $X_{18} \subset \mathbb{P} (1,1,4,6,7)$ & $\mathrm{BR}$ & 84 & $X_{36} \subset \mathbb{P} (1,7,8,9,12)$ & $\mathrm{BSR}$ \\
37 & $X_{18} \subset \mathbb{P} (1,2,3,4,9)$ & $\mathrm{BSR}$ & 85 & $X_{38} \subset \mathbb{P} (1,3,5,11,19)$ & $\mathrm{BSR}$ \\
38 & $X_{18} \subset \mathbb{P} (1,2,3,5,8)$ & $\mathrm{BR}$ & 86 & $X_{38} \subset \mathbb{P} (1,5,6,8,19)$ & $\mathrm{BSR}$ \\
39 & $X_{18} \subset \mathbb{P} (1,3,4,5,6)$ & $\mathrm{BSR}$ & 87 & $X_{40} \subset \mathbb{P} (1,5,7,8,20)$ & $\mathrm{BSR}$ \\
40 & $X_{19} \subset \mathbb{P} (1,3,4,5,7)$ & $\mathrm{BR}$ & 88 & $X_{42} \subset \mathbb{P} (1,1,6,14,21)$ & $\mathrm{BSR}$ \\
41 & $X_{20} \subset \mathbb{P} (1,1,4,5,10)$ & $\mathrm{BR}$ & 89 & $X_{42} \subset \mathbb{P} (1,2,5,14,21)$ & $\mathrm{BSR}$ \\
42 & $X_{20} \subset \mathbb{P} (1,2,3,5,10)$ & $\mathrm{BSR}$ & 90 & $X_{42} \subset \mathbb{P} (1,3,4,14,21)$ & $\mathrm{BSR}$ \\
43 & $X_{20} \subset \mathbb{P} (1,2,4,5,9)$ & $\mathrm{BR}$ & 91 & $X_{44} \subset \mathbb{P} (1,4,5,13,22)$ & $\mathrm{BSR}$ \\
44 & $X_{20} \subset \mathbb{P} (1,2,5,6,7)$ & $\mathrm{BR}$ & 92 & $X_{48} \subset \mathbb{P} (1,3,5,16,24)$ & $\mathrm{BSR}$ \\
45 & $X_{20} \subset \mathbb{P} (1,3,4,5,8)$ & $\mathrm{BR}$ & 93 & $X_{50} \subset \mathbb{P} (1,7,8,10,25)$ & $\mathrm{BSR}$ \\
46 & $X_{21} \subset \mathbb{P} (1,1,3,7,10)$ & $\mathrm{BR}$ & 94 & $X_{54} \subset \mathbb{P} (1,4,5,18,27)$ & $\mathrm{BSR}$ \\
47 & $X_{21} \subset \mathbb{P} (1,1,5,7,8)$ & $\mathrm{BR}$ & 95 & $X_{66} \subset \mathbb{P} (1,5,6,22,33)$ & $\mathrm{BSR}$ \\
48 & $X_{21} \subset \mathbb{P} (1,2,3,7,9)$ & $\mathrm{BR}$ & \\
\bottomrule
\end{tabular}
\end{table}

\begin{table}[t]
\centering
\caption{Fano threefold WCIs of codimension $2$ and index $1$}
\label{table:WCIcod2ind1} 
\begin{tabular}[]{cccccccc}
\toprule
\textnumero & $X_{d_1,d_2} \subset \mathbb{P} (a_0,\dots,a_5)$ & BP & \textnumero & $X_{d_1,d_2} \subset \mathbb{P} (a_0,\dots,a_5)$ & BP \\
\midrule
1 & $X_{2,3} \subset \mathbb{P} (1,1,1,1,1,1)$ & $\mathrm{BR}_*$ & 44 & $X_{10,12} \subset \mathbb{P} (1,2,3,5,5,7)$ & $\mathrm{BS}_*$ \\
2 & $X_{3,3} \subset \mathbb{P} (1,1,1,1,2,2)$ & $\mathrm{dP}_3$ & 45 & $X_{10,12} \subset \mathbb{P} (1,2,4,5,5,6)$ & $\mathrm{BR}$ \\
3 & $X_{3,4} \subset \mathbb{P} (1,1,1,1,2,2)$ & $\mathrm{BS}_*$ & 46 & $X_{10,12} \subset \mathbb{P} (1,3,3,4,5,7)$ & \\
4 & $X_{4,4} \subset \mathbb{P} (1,1,1,1,2,3)$ & $\mathrm{dP}_2$ & 47 & $X_{10,12} \subset \mathbb{P} (1,3,4,4,5,6)$ & $\mathrm{BSR}$ \\
5 & $X_{4,4} \subset \mathbb{P} (1,1,1,2,2,2)$ & $\mathrm{dP}_2$ & 48 & $X_{11,12} \subset \mathbb{P} (1,1,4,5,6,7)$ & $\mathrm{BS}_*$ \\
6 & $X_{4,5} \subset \mathbb{P} (1,1,1,2,2,3)$ & $\mathrm{BS}_*$ & 49 & $X_{10,14} \subset \mathbb{P} (1,1,2,5,7,9)$ & $\mathrm{BS}$ \\
7 & $X_{4,6} \subset \mathbb{P} (1,1,1,2,3,3)$ & $\mathrm{BS}_*$ & 50 & $X_{10,14} \subset \mathbb{P} (1,2,3,5,7,7)$ & $\mathrm{BS}$ \\
8 & $X_{4,6} \subset \mathbb{P} (1,1,2,2,2,3)$ & $\mathrm{BR}$ & 51 & $X_{10,14} \subset \mathbb{P} (1,2,4,5,6,7)$ & $\mathrm{BR}$ \\
9 & $X_{5,6} \subset \mathbb{P} (1,1,1,2,3,4)$ & $\mathrm{BS}_*$ & 52 & $X_{10,15} \subset \mathbb{P} (1,2,3,5,7,8)$ & $\mathrm{BS}_*$ \\
10 & $X_{5,6} \subset \mathbb{P} (1,1,2,2,3,3)$ & $\mathrm{BS}_*$ & 53 & $X_{12,13} \subset \mathbb{P} (1,3,4,5,6,7)$ & \\
11 & $X_{6,6} \subset \mathbb{P} (1,1,1,2,3,5)$ & $\mathrm{dP}_1$ & 54 & $X_{12,14} \subset \mathbb{P} (1,1,3,4,7,11)$ & \\
12 & $X_{6,6} \subset \mathbb{P} (1,1,2,2,3,4)$ & $\mathrm{dP}_1$ & 55 & $X_{12,14} \subset \mathbb{P} (1,1,4,6,7,8)$ & $\mathrm{BS}$ \\
13 & $X_{6,6} \subset \mathbb{P} (1,1,2,3,3,3)$ & $\mathrm{dP}_1$ & 56 & $X_{12,14} \subset \mathbb{P} (1,2,3,4,7,10)$ & \\
14 & $X_{6,6} \subset \mathbb{P} (1,2,2,2,3,3)$ & $\mathrm{BSR}$ & 57 & $X_{12,14} \subset \mathbb{P} (1,2,3,5,7,9)$ & $\mathrm{BS}_*$ \\
15 & $X_{6,7} \subset \mathbb{P} (1,1,2,2,3,5)$ & & 58 & $X_{12,14} \subset \mathbb{P} (1,3,4,5,7,7)$ & \\
16 & $X_{6,7} \subset \mathbb{P} (1,1,2,3,3,4)$ & $\mathrm{BS}_*$ & 59 & $X_{12,14} \subset \mathbb{P} (1,4,4,5,6,7)$ & $\mathrm{BSR}$ \\
17 & $X_{6,8} \subset \mathbb{P} (1,1,1,3,4,5)$ & $\mathrm{BS}$ & 60 & $X_{12,14} \subset \mathbb{P} (2,3,4,5,6,7)$ & $\mathrm{BSR}$ \\
18 & $X_{6,8} \subset \mathbb{P} (1,1,2,3,3,5)$ & $\mathrm{BS}_*$ & 61 & $X_{12,15} \subset \mathbb{P} (1,1,4,5,6,11)$ & $\mathrm{BS}$ \\
19 & $X_{6,8} \subset \mathbb{P} (1,1,2,3,4,4)$ & $\mathrm{BS}$ & 62 & $X_{12,15} \subset \mathbb{P} (1,3,4,5,6,9)$ & $\mathrm{BS}$ \\
20 & $X_{6,8} \subset \mathbb{P} (1,2,2,3,3,4)$ & $\mathrm{BR}$ & 63 & $X_{12,15} \subset \mathbb{P} (1,3,4,5,7,8)$ & $\mathrm{BS}_*$ \\
21 & $X_{6,9} \subset \mathbb{P} (1,1,2,3,4,5)$ & $\mathrm{BS}_*$ & 64 & $X_{12,16} \subset \mathbb{P} (1,2,5,6,7,8)$ & $\mathrm{BR}$ \\
22 & $X_{7,8} \subset \mathbb{P} (1,1,2,3,4,5)$ & $\mathrm{BS}_*$ & 65 & $X_{14,15} \subset \mathbb{P} (1,2,3,5,7,12)$ & \\
23 & $X_{6,10} \subset \mathbb{P} (1,1,2,3,5,5)$ & $\mathrm{BS}$ & 66 & $X_{14,15} \subset \mathbb{P} (1,2,5,6,7,9)$ & \\
24 & $X_{6,10} \subset \mathbb{P} (1,2,2,3,4,5)$ & $\mathrm{BR}$ & 67 & $X_{14,15} \subset \mathbb{P} (1,3,4,5,7,10)$ & \\
25 & $X_{8,9} \subset \mathbb{P} (1,1,2,3,4,7)$ & & 68 & $X_{14,15} \subset \mathbb{P} (1,3,5,6,7,8)$ & \\
26 & $X_{8,9} \subset \mathbb{P} (1,1,3,4,4,5)$ & $\mathrm{BS}_*$ & 69 & $X_{14,16} \subset \mathbb{P} (1,1,5,7,8,9)$ & $\mathrm{BS}$ \\
27 & $X_{8,9} \subset \mathbb{P} (1,2,3,3,4,5)$ & & 70 & $X_{14,16} \subset \mathbb{P} (1,3,4,5,7,11)$ & \\
28 & $X_{8,10} \subset \mathbb{P} (1,1,2,3,5,7)$ & $\mathrm{BS}_*$ & 71 & $X_{14,16} \subset \mathbb{P} (1,4,5,6,7,8)$ & $\mathrm{BSR}$ \\
29 & $X_{8,10} \subset \mathbb{P} (1,1,2,4,5,6)$ & $\mathrm{BS}$ & 72 & $X_{15,16} \subset \mathbb{P} (1,2,3,5,8,13)$ & \\
30 & $X_{8,10} \subset \mathbb{P} (1,1,3,4,5,5)$ & $\mathrm{BS}$ & 73 & $X_{15,16} \subset \mathbb{P} (1,3,4,5,8,11)$ & \\
31 & $X_{8,10} \subset \mathbb{P} (1,2,3,4,4,5)$ & $\mathrm{BR}$ & 74 & $X_{14,18} \subset \mathbb{P} (1,2,3,7,9,11)$ & $\mathrm{BS}$ \\
32 & $X_{9,10} \subset \mathbb{P} (1,1,2,3,5,8)$ & & 75 & $X_{14,18} \subset \mathbb{P} (1,2,6,7,8,9)$ & $\mathrm{BR}$ \\
33 & $X_{9,10} \subset \mathbb{P} (1,1,3,4,5,6)$ & $\mathrm{BS}_*$ & 76 & $X_{12,20} \subset \mathbb{P} (1,4,5,6,7,10)$ & $\mathrm{BR}$ \\
34 & $X_{9,10} \subset \mathbb{P} (1,2,2,3,5,7)$ & & 77 & $X_{16,18} \subset \mathbb{P} (1,1,6,8,9,10)$ & $\mathrm{BS}$ \\
35 & $X_{9,10} \subset \mathbb{P} (1,2,3,4,5,5)$ & & 78 & $X_{16,18} \subset \mathbb{P} (1,4,6,7,8,9)$ & $\mathrm{BSR}$ \\
36 & $X_{8,12} \subset \mathbb{P} (1,1,3,4,5,7)$ & $\mathrm{BS}$ & 79 & $X_{18,20} \subset \mathbb{P} (1,4,5,6,9,14)$ & \\
37 & $X_{8,12} \subset \mathbb{P} (1,2,3,4,5,6)$ & $\mathrm{BR}$ & 80 & $X_{18,20} \subset \mathbb{P} (1,4,5,7,9,13)$ & \\
38 & $X_{9,12} \subset \mathbb{P} (1,2,3,4,5,7)$ & $\mathrm{BS}_*$ & 81 & $X_{18,20} \subset \mathbb{P} (1,5,6,7,9,11)$ & \\
39 & $X_{10,11} \subset \mathbb{P} (1,2,3,4,5,7)$ & & 82 & $X_{18,22} \subset \mathbb{P} (1,2,5,9,11,13)$ & $\mathrm{BS}$ \\
40 & $X_{10,12} \subset \mathbb{P} (1,1,3,4,5,9)$ & & 83 & $X_{20,21} \subset \mathbb{P} (1,3,4,7,10,17)$ \\
41 & $X_{10,12} \subset \mathbb{P} (1,1,3,5,6,7)$ & $\mathrm{BS}$ & 84 & $X_{18,30} \subset \mathbb{P} (1,6,8,9,10,15)$ & $\mathrm{BSR}$ \\
42 & $X_{10,12} \subset \mathbb{P} (1,1,4,5,6,6)$ & $\mathrm{BS}$ & 85 & $X_{24,30} \subset \mathbb{P} (1,8,9,10,12,15)$ & $\mathrm{BSR}$ \\
43 & $X_{10,12} \subset \mathbb{P} (1,2,3,4,5,8)$ \\
\bottomrule
\end{tabular}
\end{table}

\begin{table}[t]
\centering
\caption{Fano threefold weighted hypersurfaces of index $> 1$}
\label{table:WCIcod1indhigh} 
\begin{tabular}[]{cccccccccc}
\toprule
\textnumero & $X_d \subset \mathbb{P} (a_0,\dots,a_4)$ & BS & Rat & \textnumero & $X_d \subset \mathbb{P} (a_0,\dots,a_4)$ & BS & Rat \\
\midrule
96 & $X_3 \subset \mathbb{P} (1,1,1,1,1)$ & $-$ & $-$ & 114 & $X_6 \subset \mathbb{P} (1,1,2,3,4)$ & $-$ & $+$ \\
97 & $X_4 \subset \mathbb{P} (1,1,1,1,2)$ & $-$ & $-$ & 115 & $X_6 \subset \mathbb{P} (1,2,2,3,3)$ & $-$ & $+$ \\
98 & $X_{6} \subset \mathbb{P} (1,1,1,2,3)$ & $-$ & $-$ & 116 & $X_{10} \subset \mathbb{P} (1,2,3,4,5)$ & $-$ & $-$ \\
99 & $X_{10} \subset \mathbb{P} (1,1,2,3,5)$ & $-$ & ? & 117 & $X_{15} \subset \mathbb{P} (1,3,4,5,7)$ & $-$ & ? \\
100 & $X_{18} \subset \mathbb{P} (1,2,3,5,9)$ & $+$ & $-$ & 118 & $X_6 \subset \mathbb{P} (1,1,2,3,5)$ & $-$ & $+$ \\
101 & $X_{22} \subset \mathbb{P} (1,2,3,7,11)$ & $+$ & $-$ & 119 & $X_6 \subset \mathbb{P} (1,2,2,3,5)$ & $-$ & $+$ \\
102 & $X_{26} \subset \mathbb{P} (1,2,5,7,13)$ & $+$ & $-$ & 120 & $X_6 \subset \mathbb{P} (1,2,3,3,4)$ & $-$ & $+$ \\
103 & $X_{38} \subset \mathbb{P} (2,3,5,11,19)$ & $+$ & $-$ & 121 & $X_8 \subset \mathbb{P} (1,2,3,4,5)$ & $-$ & $+$ \\
104 & $X_2 \subset \mathbb{P} (1,1,1,1,1)$ & $-$ & $+$ & 122 & $X_{14} \subset \mathbb{P} (2,3,4,5,7)$ & $-$ & ? \\
105 & $X_3 \subset \mathbb{P} (1,1,1,1,2)$ & $-$ & $+$ & 123 & $X_6 \subset \mathbb{P} (1,2,3,3,5)$ & $-$ & $+$ \\
106 & $X_4 \subset \mathbb{P} (1,1,1,2,2)$ & $-$ & $+$ & 124 & $X_{10} \subset \mathbb{P} (1,2,3,5,7)$ & $-$ & $+$ \\
107 & $X_6 \subset \mathbb{P} (1,1,2,2,3)$ & $-$ & $-$ & 125 & $X_{12} \subset \mathbb{P} (1,3,4,5,7)$ & $-$ & $+$ \\
108 & $X_{12} \subset \mathbb{P} (1,2,3,4,5)$ & $-$ & ? & 126 & $X_6 \subset \mathbb{P} (1,2,3,4,5)$ & $-$ & $+$ \\
109 & $X_{15} \subset \mathbb{P} (1,2,3,5,7)$ & $-$ & ? & 127 & $X_{12} \subset \mathbb{P} (2,3,4,5,7)$ & $-$ & $+$ \\
110 & $X_{21} \subset \mathbb{P} (1,3,5,7,8)$ & $+$ & $-$ & 128 & $X_{12} \subset \mathbb{P} (1,4,5,6,7)$ & $-$ & $+$ \\
111 & $X_4 \subset \mathbb{P} (1,1,1,2,3)$ & $-$ & $+$ & 129 & $X_{10} \subset \mathbb{P} (2,3,4,5,7)$ & $-$ & $+$ \\
112 & $X_6 \subset \mathbb{P} (1,1,2,3,3)$ & $-$ & $+$ & 130 & $X_{12} \subset \mathbb{P} (3,4,5,6,7)$ & $-$ & $+$ \\
113 & $X_4 \subset \mathbb{P} (1,1,2,2,3)$ & $-$ & $+$ & \\
\bottomrule
\end{tabular}
\end{table}

\FloatBarrier
\begingroup
\renewcommand\arraystretch{1.2}
\begin{longtable}{rccl}
\caption{Fano threefold WCIs of codimension $2$ and index $> 1$}
\label{table:WCIcod2highind} \\

\toprule
\textnumero & $X_{d_1,d_2} \subset \mathbb{P} (a_0,\dots,a_5)$ & Centre &Models \\
\midrule
\endfirsthead

\toprule
\textnumero & $X_{d_1,d_2} \subset \mathbb{P} (a_0,\dots,a_5)$ & Centre &Models \\
\midrule
\endhead

\midrule
\multicolumn{4}{r}{Continued on the next page.} \\ 
\endfoot

\bottomrule
\endlastfoot

\multirow{1}{*}{87} & \multirow{1}{*}{$X_{4,4} \subset \mathbb{P} (1,1,1,2,2,3)$} & $\frac{1}{3}(1,1,2)$ & $Y/\mathbb P^2$ \\

\multirow{1}{*}{88} &  \multirow{1}{*}{$X_{4,6} \subset \mathbb{P} (1,1,2,2,3,3)$} & $2\times \frac{1}{3}(1,1,2)$ & $\mathrm{dP}_3$ \\

\multirow{1}{*}{89} &  \multirow{1}{*}{$X_{6,6} \subset \mathbb{P} (1,1,2,2,3,5)$} & $\frac{1}{5}(1,1,4)$ & $\mathrm{dP}_3$ \\

\multirow{1}{*}{90} & \multirow{1}{*}{$X_{6,8} \subset \mathbb{P} (1,1,2,3,4,5)$} &$\frac{1}{5}(1,2,3)$ & $\mathrm{dP}_2$ \\

\multirow{1}{*}{91} & \multirow{1}{*}{$X_{6,6} \subset \mathbb{P} (1,2,2,3,3,3)$} &$4\times \frac{1}{3}(1,1,2)$ & $X_{3,3}\subset \mathbb P(1,1,1,1,1,2)$ \\

\multirow{1}{*}{92} & \multirow{1}{*}{$X_{6,8} \subset \mathbb{P} (1,2,2,3,3,5)$} & $\frac{1}{5}(1,1,4)$ & $X_4\subset \mathbb P^4$ \\ & & $2\times \frac{1}{3}(1,1,2)$ & $X_{3,4} \subset \mathbb P(1,1,1,1,2,2)$ \\

\multirow{1}{*}{93} & \multirow{1}{*}{$X_{6,10} \subset \mathbb{P} (1,2,2,3,5,5)$} & $2\times \frac{1}{5}(1,1,4)$ & $X_5\subset \mathbb P^(1,1,1,1,2)$  \\

\multirow{1}{*}{94} & \multirow{1}{*}{$X_{8,10} \subset \mathbb{P} (1,2,2,3,3,5)$} & $\frac{1}{7}(1,1,6)$ & $ X_5\subset \mathbb P(1,1,1,1,2)$ \\ & & $\frac{1}{3}(1,1,2)$ & $ X_{4,5} \subset \mathbb P(1,1,1,2,2,3)$ \\
     
\multirow{1}{*}{95} & \multirow{1}{*}{$X_{8,10} \subset \mathbb{P} (1,2,3,4,5,5)$} & $2\times \frac{1}{5}(1,2,3)$ & $X_{4,5}\subset \mathbb P(1,1,1,2,2,3)$ \\ & & $\frac{1}{3}(1,1,2)$ & $X_{4,6} \subset \mathbb P(1,1,2,2,2,3)$ \\

\multirow{1}{*}{96} & \multirow{1}{*}{$X_{8,12} \subset \mathbb{P} (1,2,3,4,5,7)$} & $\frac{1}{7}(1,2,5)$ & $X_{6}\subset \mathbb P(1,1,1,2,2)$ \\ & & $\frac{1}{5}(1,2,3)$ & $X_{4,6} \subset \mathbb P(1,1,1,2,3,3)$ \\

\multirow{1}{*}{97} & \multirow{1}{*}{$X_{10,14} \subset \mathbb{P} (1,2,2,5,7,9)$} & $\frac{1}{9}(1,1,8)$ & $X_{7}\subset \mathbb P(1,1,1,2,3)$ \\ 

\multirow{1}{*}{98} & \multirow{1}{*}{$X_{10,12} \subset \mathbb{P} (1,2,3,4,5,9)$} & $\frac{1}{9}(1,2,7)$ & $X_{6}\subset \mathbb P(1,1,1,2,2)$ \\ & & $\frac{1}{3}(1,1,2)$ & $X_{6,6} \subset \mathbb P(1,1,1,2,3,3)$ \\

\multirow{1}{*}{99} & \multirow{1}{*}{$X_{10,12} \subset \mathbb{P} (1,2,3,5,6,7)$} & $\frac{1}{7}(1,3,4)$ & $X_{5,6}\subset \mathbb P(1,1,1,2,3,4)$ \\ & & $2\times \frac{1}{3}(1,1,2)$ & $X \subset \mathbb P(1,1,2,3,3,4,7,10)$ \\

\multirow{1}{*}{100} & \multirow{1}{*}{$X_{12,14} \subset \mathbb{P} (1,2,3,4,7,11)$} & $\frac{1}{11}(1,2,9)$ & $X_{7}\subset \mathbb P(1,1,1,2,3)$ \\ 

\multirow{1}{*}{101} & \multirow{1}{*}{$X_{10,12} \subset \mathbb{P} (2,2,3,5,5,7)$} & $\frac{1}{7}(1,1,6)$ & $X_{6}\subset \mathbb P(1,1,1,1,3)$ \\ & & $2\times \frac{1}{5}(1,1,4)$ & $X_{5,6} \subset \mathbb P(1,1,1,2,3,4)$ \\
     
\multirow{1}{*}{102} & \multirow{1}{*}{$X_{10,14} \subset \mathbb{P} (2,2,3,5,7,7)$} & $2\times \frac{1}{7}(1,1,6)$ & $X_{7}\subset \mathbb P(1,1,1,2,3)$ \\ & & $\frac{1}{3}(1,1,2)$ & $X_{12} \subset \mathbb P(1,1,1,4,6)$ \\

\multirow{1}{*}{103} & \multirow{1}{*}{$X_{10,12} \subset \mathbb{P} (2,3,3,4,5,7)$} & $\frac{1}{7}(1,2,5)$ & $\mathrm{dP}_1$ \\ & & $4\times \frac{1}{3}(1,1,2)$ & $X \subset \mathbb P(1,2,3,3,4,5,7,9)$ \\

\multirow{1}{*}{104} & \multirow{1}{*}{$X_{14,16} \subset \mathbb{P} (1,2,5,7,8,9)$} & $\frac{1}{9}(1,4,5)$ & $X_{7,8} \subset \mathbb P(1,1,2,3,4,5)$ \\ & & $\frac{1}{5}(1,1,4)$ & $X \subset \mathbb P(1,1,3,4,4,5,9,13)$ \\

\multirow{1}{*}{105} & \multirow{1}{*}{$X_{12,14} \subset \mathbb{P} (2,2,3,5,7,9)$} & $\frac{1}{9}(1,1,8)$ & $X_{7} \subset \mathbb P(1,1,1,2,3)$ \\ & & $\frac{1}{5}(1,1,4)$ & $X_{6,7} \subset \mathbb P(1,1,2,3,3,4)$ \\
& & $\frac{1}{3}(1,1,2)$ & $X_{12} \subset \mathbb P(1,1,2,4,5)$ \\

\multirow{1}{*}{106} & \multirow{1}{*}{$X_{14,18} \subset \mathbb{P} (2,2,3,7,9,11)$} & $\frac{1}{11}(1,1,10)$ & $X_{9} \subset \mathbb P(1,1,2,3,3)$ \\ & & $2\times \frac{1}{3}(1,1,2)$ & $X \subset \mathbb P(1,1,6,8,9,10,11,12)$ \\

\multirow{1}{*}{107} & \multirow{1}{*}{$X_{12,14} \subset \mathbb{P} (2,3,4,5,7,7)$} & $\frac{1}{7}(1,2,5)$ & $X_{6,7} \subset \mathbb P(1,1,2,2,3,5)$ \\ & & $\frac{1}{5}(1,1,4)$ & $X_{6,10} \subset \mathbb P(1,2,2,3,4,5)$ \\

\multirow{1}{*}{108} & \multirow{1}{*}{$X_{14,16} \subset \mathbb{P} (1,2,3,4,7,11)$} & $\frac{1}{11}(1,2,9)$ & $X_{8}\subset \mathbb P(1,1,2,2,3)$ \\ 

\multirow{1}{*}{109} & \multirow{1}{*}{$X_{18,22} \subset \mathbb{P} (2,2,5,9,11,13)$} & $\frac{1}{13}(1,1,12)$ & $ X_{11}\subset \mathbb P(1,1,2,3,5)$ \\\

\multirow{1}{*}{110} & \multirow{1}{*}{$X_{18,20} \subset \mathbb{P} (2,4,5,7,9,13)$} & $\frac{1}{13}(1,2,11)$ & $X_{10}\subset \mathbb P(1,1,2,2,5)$ \\

\multirow{1}{*}{111} & \multirow{1}{*}{$X_{18,20} \subset \mathbb{P} (2,5,6,7,9,11)$} & $\frac{1}{11}(1,3,8)$ & $X_{9,10} \subset \mathbb P(1,1,2,3,5,8)$ \\ & & $\frac{1}{7}(1,1,6)$ & $X \subset \mathbb P(1,3,5,6,7,8,11,14)$ \\

\multirow{1}{*}{112} & \multirow{1}{*}{$X_{6,6} \subset \mathbb{P} (1,1,2,3,3,5)$} & $\frac{1}{5}(1,1,4)$ & $Y/\mathbb P(1,1,2)$ \\

\multirow{1}{*}{113} & \multirow{1}{*}{$X_{6,6} \subset \mathbb{P} (1,2,2,3,3,4)$} & $\frac{1}{4}(1,1,3)$ & $X_{4,4}\subset \mathbb P(1,1,1,2,2,3)$ \\ & & $\frac{1}{2}(1,1,1)$ & $\mathrm{dP}_4$ \\

\multirow{1}{*}{114} & \multirow{1}{*}{$X_{6,9} \subset \mathbb{P} (1,2,3,3,4,5)$} & $\frac{1}{5}(1,1,4)$ & $\mathrm{dP}_3$ \\ & & $\frac{1}{4}(1,1,3)$ & $X_{4,6} \subset \mathbb P(1,1,2,2,3,3)$ \\
& & $\frac{1}{2}(1,1,1)$ & $ X_{2,3} \subset \mathbb P^5$ \\

\multirow{1}{*}{115} & \multirow{1}{*}{$X_{12,15} \subset \mathbb{P} (1,3,4,5,6,11)$} & $\frac{1}{11}(1,2,9)$ & $X_5 \subset \mathbb P(1,1,1,1,2)$ \\

\multirow{1}{*}{116} & \multirow{1}{*}{$X_{9,12} \subset \mathbb{P} (2,3,3,4,5,7)$} & $\frac{1}{7}(1,1,6)$ & $\mathrm{dP}_2$ \\ & & $\frac{1}{5}(1,1,4)$ & $X_{3,4} \subset \mathbb P(1,1,1,2,2,3)/\bm{\mu}_2$ \\
& & $\frac{1}{2}(1,1,1)$ & $X_{6} \subset \mathbb P(1,1,1,2,2)$ \\

\multirow{1}{*}{117} & \multirow{1}{*}{$X_{12,15} \subset \mathbb{P} (3,3,4,5,7,8)$} & $\frac{1}{8}(1,1,7)$ & $X_{8,10}\subset \mathbb P(1,2,2,3,5,7)$ \\ & & $\frac{1}{7}(1,1,6)$ & $X_{3,4} \subset \mathbb P(1,1,1,2,2,3)/\bm{\mu}_2$ \\
& & $\frac{1}{4}(1,1,3)$ & $X_{8,10} \subset \mathbb P(1,2,2,3,5,7)$ \\

\multirow{1}{*}{118} & \multirow{1}{*}{$X_{6,8} \subset \mathbb{P} (1,2,3,3,4,5)$} & $\frac{1}{5}(1,2,3)$ & $X_6 \subset \subset \mathbb P(1,1,2,2,3)$ \\ & & $\frac{1}{3}(1,1,2)$ & $Y/\mathbb P(1,2,3)$ \\

\multirow{1}{*}{119} & \multirow{1}{*}{$X_{8,10} \subset \mathbb{P} (1,2,3,4,5,7)$} & $\frac{1}{7}(1,3,4)$ & $Y/\mathbb P(1,2,3)$ \\ & & $\frac{1}{3}(1,1,2)$ & $\mathrm{dP}_4$ \\

\multirow{1}{*}{120} & \multirow{1}{*}{$X_{8,12} \subset \mathbb{P} (1,3,4,4,5,7)$} & $\frac{1}{7}(1,1,6)$ & $\mathrm{dP}_3$ \\ & & $\frac{1}{5}(1,1,4)$ & $ X_{6,9}\subset \mathbb P(1,2,3,3,4,5)$ \\

\multirow{1}{*}{121} & \multirow{1}{*}{$X_{10,12} \subset \mathbb{P} (1,3,4,5,6,7)$} & $\frac{1}{7}(1,2,5)$ & $\mathrm{dP}_3$ \\ & & $2\times \frac{1}{3}(1,1,2)$ & $X_{2,3}\subset \mathbb P^5$ \\

\multirow{1}{*}{122} & \multirow{1}{*}{$X_{10,12} \subset \mathbb{P} (2,3,4,5,5,7)$} & $\frac{1}{7}(1,3,4)$ & $\mathrm{dP}_3$ \\ & & $\frac{1}{5}(1,2,3)$ & $X_{4,6}\subset \mathbb P(1,1,2,2,3,3)$ \\

\multirow{1}{*}{123} & \multirow{1}{*}{$X_{12,14} \subset \mathbb{P} (2,3,4,5,7,9)$} & $\frac{1}{9}(1,4,5)$ & $X_{6,6}\subset \mathbb P(1,1,2,2,3,5)$ \\ & & $\frac{1}{5}(1,2,3)$ & $ X_{6,6} \subset \mathbb P(1,2,2,3,3,3)$ \\
& & $\frac{1}{3}(1,1,2)$ & $X_{3,4} \subset \mathbb P(1,1,1,1,2,2)/\bm{\mu}_2$ \\

\multirow{1}{*}{124} & \multirow{1}{*}{$X_{18,20} \subset \mathbb{P} (4,5,6,7,9,11)$} & $\frac{1}{11}(1,4,7)$ & $X_{15}\subset \mathbb P(1,2,3,5,7)$ \\ & & $\frac{1}{7}(1,3,4)$ & $X_{10} \subset \mathbb P(1,1,2,3,5)$ \\

\multirow{1}{*}{125} & \multirow{1}{*}{$X_{10,15} \subset \mathbb{P} (2,3,5,5,7,8)$} & $\frac{1}{8}(1,1,7)$ & $X_{8,12}\subset \mathbb P(1,3,4,4,5,7)$ \\ & & $\frac{1}{7}(1,1,6)$ & $X_{4,6} \subset \mathbb P(1,1,2,2,3,3)/\bm{\mu}_2$ \\
& & $\frac{1}{2}(1,1,1)$ & $X_{8,12} \subset \mathbb P(1,3,4,4,5,7)$ \\
\end{longtable}
\endgroup
\FloatBarrier

\bibliography{literatur}
\bibliographystyle{alpha}

\newpage

\end{document}